\documentclass[12pt]{article}
\usepackage{geometry}
\usepackage{amssymb}
\usepackage{color}

\usepackage{amsmath,amsthm}
\allowdisplaybreaks

\newif\ifleft
\lefttrue 

\newif\ifpre
\prefalse

\newif\ifhide
\hidefalse

\def\1{\mathbf{1}}

\def\P{\mathbb{P}}

\newtheorem{lemma} {Lemma}
\newtheorem{theorem} {Theorem}
\newtheorem{Proposition}{Proposition}
\newtheorem{Definition}{Definition}
\newtheorem{remark} {Remark}

\newif \ifshowup
\showupfalse 

\title{{\normalsize\tt\hfill\jobname.tex}\\
      \bf Existence and
uniqueness theorems for solutions of McKean--Vlasov
stochastic equations \\
{\small In memory of A.V. Skorokhod (10.09.1930 --
03.01.2011)}\\~\\
Yu. S. Mishura\footnote{Taras Shevchenko National University of Kyiv; myus@univ.kiev.ua}, \, \& \, A. Yu. Veretennikov\footnote{University of Leeds, UK
$\,$ \& $\,$ National Research University Higher School of Economics, Russian Federation, and Institute for Information Transmission Problems, Moscow, Russia; a.veretennikov@leeds.ac.uk}
}

\begin{document}
\maketitle

\begin{abstract}
New weak and strong existence and weak and strong uniqueness results for multi-dimensional stochastic McKean--Vlasov equation are established  under relaxed regularity conditions. Weak existence requires a non-degeneracy of diffusion and no more than a linear growth of both coefficients in the state variable. Weak and strong uniqueness are established under the restricted assumption of diffusion, yet without any regularity of the drift; this part is based on the analysis of the total variation metric.
\end{abstract}

\maketitle

\section{Introduction}
 Our subject is solutions  of the stochastic It\^o-McKean-Vlasov equation in $\mathbb R^d$
\begin{equation}\label{e1}
dX_t = B[t,X_t, \mu_t]dt + \Sigma[t,X_t, \mu_t]dW_t, \;\; t\ge 0, \qquad
X_0=x_0,
\end{equation}
in a particular situation called ``true McKean-Vlasov case''  under the convention
\begin{equation}\label{e200}
B[t,x,\mu]=\int
b(t,x,y)\mu(dy), \;\; \Sigma[t,x,\mu]=\int
\sigma(t,x,y)\mu(dy),
\end{equation}
and under certain non-degeneracy assumptions.
Here $W$ is a standard $d_1$-dimensional
Wiener process, $b$ and $\sigma$ are vector and matrix Borel
functions of corresponding dimensions $d$ and $d\times
d_1$, $\mu_t$ is the distribution of the process $X$ at time
$t$. The initial data $x_0$ may be random  but independent of
$W$; a non-random value is also allowed. Historically,  Vlasov's idea, proposed originally in 1938 and contained in    the reprinted paper \cite{Vlasov68}, called mean field interaction in
mathematical physics and stochastic analysis, assumes
that
for a large multi-particle ensemble with ``weak interaction''
between particles, this interaction for one particle
with others may be effectively replaced by an averaged
field. A class of  equations  of type (\ref{e1}) was proposed by M.~Kac
\cite{Kac}	
as a stochastic ``toy model'' for the Vlasov kinetic equation of
plasma. The systematic study of such equations
was started by McKean \cite{McKean}. The reference \cite{Sz} provides an introduction to the
whole area with links to the paper \cite{Dobr} as the most important preceding background deterministic
paper.

McKean--Vlasov's equations, being clearly more
involved than It\^o's SDEs, arise in
multi-agent systems (see \cite{Ben, Bossy_Talay}),  as well as in some other areas of high interest
such as filtering (see \cite{Crisan_Xiong}). These processes also very closely relate to so called self-stabilizing processes (diffusions, in particular), which is, actually, another name for non-linear diffusions in the ``ergodic'' situation, (see \cite{Herrmann_Tugaut}). In what concerns ``propagation of chaos''
for the equation (\ref{e1}), we refer the reader to
\cite{Sz} and \cite[Theorem 4.3]{Chiang}.  In the authors' view, it may be fruitful to separate different
aspects, including time
discretization and ``propagation of chaos'' for
multi-particle case, and to consider
approximations differently from the basic existence and
uniqueness issues; only the latter two are the main subjects of the present paper.
Many control problems
lead to discontinuous coefficients. This is one of motivations for looking for
existence and uniqueness under minimal regularity of the coefficients.

As to earlier works in this area, one of the most important papers  is \cite{Funaki} where the martingale problem for a similar McKean-Vlasov SDE is tackled. It is not very easy to compare our regularity assumptions with those in \cite{Funaki} because the latter are given not directly  in terms of coefficients (please compare with   (2.1) in the Assumption I from \cite{Funaki}). We do not assume continuity with respect to the state variable $x$ replacing it by the non-degeneracy of the diffusion matrix. Neither our linear growth bound is comparable directly with the Lyapunov type conditions in \cite{Funaki}.
More
general growth conditions were studied in \cite{Chiang};
however, our regularity conditions admit
just measurable coefficients in $x$, especially, for weak existence, and, hence,
overall, our results are not
covered by \cite{Chiang} either.

Our goal is to establish weak existence analogues to Krylov's
weak existence for It\^o's equations which is more general than in
earlier papers. A more general equation  is tackled with a
possibly non-square matrix $\sigma$, which may be useful in
applications and which case was not covered in
\cite{Chiang}. Further, we propose a different method which could be of interest in some other settings. In the homogeneous case and under less general conditions, using a different technique, weak existence and weak uniqueness were established in  \cite{BJ97} and \cite{BJSM98}. In \cite{Ver06} there is a result on strong existence for the equation similar to (\ref{e1}) only with a unit matrix diffusion; however, strong and weak uniqueness, along with ``propagation of chaos'', i.e., with  convergence of particle approximations, are established there under restrictive additional assumptions on the drift which include Lipschitz and some other conditions. In the present paper, weak
and strong uniqueness are established for linear growing in the state variable and
measurable drifts under additional assumptions on the (variable) diffusion coefficient.

In applications  where some additional regularization by white noise is often required it may be useful to have a result for references with dimensions $d_1 \ge d$ rather than just for $d_1=d$. This case is rarely tackled in the literature and it is not easy to find a suitable reference; this was the main reason why we included this extension. Despite the widespread intuitive belief that for weak solutions or weak uniqueness everything which may be desired only depends on the matrix $\sigma\sigma^*$, in fact, conditions in the McKean--Vlasov case usually do require certain properties of $\sigma$, not $\sigma\sigma^*$ (please compare, for example, with \cite{Funaki}). Hence, even if some results for $d_1 = d$ can be extended to $d_1\ge d$, yet it is not automatic. Unlike the setting in the paper \cite{Chiang}, we allow non-homogeneous
coefficients depending on time; a formal
reduction to a homogeneous case by considering
a couple $(t,X_t)$ would require unnecessary additional
conditions due to the degeneracy. Our   method of proof  is also different
from that used in \cite{Chiang}:
we use explicitly Skorokhod's single probability space
approach combined with Krylov's integral estimates for
It\^o's processes.

Strong existence for McKean--Vlasov equations in our paper is derived from strong
existence for ``ordinary'' or ``linearized'' It\^o's equations with a fixed flow $(\mu_t)$. The famous Yamada--Watanabe principle (see \cite{IkedaWatanabe}, \cite{Kurtz},  \cite{YamadaWatanabe}, \cite{KZ}) concerning weak existence and pathwise uniqueness here has a remote analogue in terms of the equivalence of weak and strong uniqueness, yet, under additional assumptions.
In all main results of the paper (but not in Propositions \ref{pro22} and \ref{thm4}) it is assumed that the drift, and in the Theorem \ref{thm1} diffusion as well, satisfies a
linear growth bound condition.
The linear growth is useful
because of numerous applications where,
at least, the drift is often unbounded; further extensions on a faster
non-linear growth
usually require Lyapunov type conditions, which are not
considered in this paper.

The structure of the paper is as follows.   In the Section \ref{sec:we}  weak existence is established.    Theorem \ref{thm1}  there
mimics Krylov's weak existence result for It\^o's SDEs from \cite{Kry69} for a homogeneous case, and from \cite{Kry} for a  non-homogeneous case; see also \cite{K-V}. No regularity  of the coefficients is assumed with respect to the state variable $x$. The proof is split into three parts. The first two parts, given in the Proposition \ref{pro22} and Proposition \ref{thm4}, are devoted to the  case  under a bit restrictive additional assumptions on the diffusion; the third part extends the consideration to the general situation, i.e. to a not necessarily quadratic and symmetric diffusion matrix. Section \ref{sec:se} is devoted to strong solutions and to weak and strong uniqueness.  Weak
uniqueness and strong uniqueness are established simultaneously  under identical (for weak
and for strong uniqueness) sets of conditions.
The latter do involve some restriction on the diffusion
coefficient which should  not depend on the measure
in  the Theorem \ref{thm5}.
For a completeness of the paper, two important classical Skorokhod     lemmas (Lemma \ref{app2} and Lemma \ref{app1}) are provided in   Appendix (Section \ref{Appe1}), along with a  ``localized'' version of certain Krylov  bounds (Lemma \ref{lekrybd}).
We shall use the following abbreviations for inequalities: by CBS we encrypt Cauchy-Buniakovsky-Schwarz inequality    and BCM denotes Bienaym\'e-Chebyshev-Markov inequality.

\section{Weak existence}\label{sec:we}
\subsection{Main results}
Before we turn to the main results, let us recall the definitions and a fact from functional analysis.
 \begin{Definition}\label{Def1}
The  triple $(X_t, \mu_t, W_t)$ is called solution of the equation (\ref{e1}) iff $(W_t)$ is a $d_1$-dimensional Wiener process  with a filtration $({\cal F}_t)$ such that for each $t$,  $X_t$ is ${\cal F}_t$-measurable, $X_t$ is continuous in $t$,  and
$$
\mathbb P\left(
X_t - x_0 - \int_0^t B[s,X_s, \mu_s]ds - \int_0^t  \Sigma[s,X_s, \mu_s]dW_s = 0, \; t\ge 0
\right)=1,
$$
with all the integrals   under the sign of probability being  correctly defined, and with $\mu_t$ being a marginal distribution of $X_t$ for each $t\ge 0$.
This solution is called strong iff for each $t$ the random variable $X_t$ is measurable with respect to the sigma-algebra ${\cal F}^W_t$ (sigma-algebra generated by Wiener process $W$); all other solutions are called weak.
\end{Definition}
Note that in the case of strong solution, it exists on any probability space with a $d_1$-dimensional Wiener process $W$. Following a tradition of It\^o SDE theory and slightly abusing a rigorous wording in the definition above, we will usually call solution just the first component $X_t$ of the triple $(X_t, \mu_t, W_t)$ yet with a compulsory property that $\mu_t$ is a marginal distribution of $X_t$ for each $t$. The next lemma is probably a common knowledge since it is never mentioned in the papers on the topic. Yet, it seems that this result does require some integrability conditions; so it is stated here  with a brief sketch of the proof so as to make sure that the linear growth conditions suffice.
Let $|\cdot|$ stand  for the Euclidean norm for any vector in $\mathbb{R}^d$,  and $\|\cdot \|$ for the standard  matrix norm, namely, $\|\sigma\| =
\left(\sum_{i,j}\sigma_{ij}^2\right)^{1/2}$.

\begin{lemma}\label{cor1}
In terms of notation  (\ref{e200}), let the Borel coefficients $b(t,x,y)$ and $\sigma(t,x,y)$ for each $x$ satisfy $$\sup_{t,y}(|b(t,x,y)| + \|\sigma(t,x,y)\|)\le C(x)$$ with some locally bounded Borel function $C(x), \, x\in \mathbb R^d$, and let $\mu_t(dy)$ be a marginal distribution of any solution $X_t$ of the equation (\ref{e1}).  Then the functions $\widetilde b(t,x):=
B[t, x, \mu_t]$ and $\widetilde \sigma(t,x):=
\Sigma[t, x, \mu_t]$ are Borel measurable in $(t, x)$.
\end{lemma}

\noindent
\begin{proof} Let $(X_t, \mu_t, W_t)$ be a solution of (\ref{e1}) on some probability space $(\Omega, {\cal F}, \mathbb P)$ with a $d_1$-dimensional Wiener process $W$, and consider another {\em independent} solution $(\xi_t, \mu_t, W_t')$ {\em with the same marginal distribution $\mu_t$ of $\xi_t$}, say, on another probability space $(\Omega', {\cal F}', \mathbb P')$ with a $d_1$-dimensional Wiener process $W'$:
\begin{equation}\label{exi}
d\xi_t = B[t,\xi_t, \mu_t]dt + \Sigma[t,\xi_t, \mu_t]dW'_t, \;\; t\ge 0, \qquad
{\cal L}(\xi_0)={\cal L}(x_0).
\end{equation}
Then the coefficient $B[t,x,\mu_t]$ can be written as
\[
B[t,x,\mu_t] = \mathbb E' b(t,x,\xi_t),
\]
where $\mathbb E'$ stands for expectation with respect to the probability measure $\mathbb P'$. (Later on we will be using the  notation $\mathbb E^3$ instead of $\mathbb E'$.)
Now, the function $b(t,x,y)$ is Borel measurable in $(t,x,y)$ by the assumption, and the function $\xi_t(\omega')$ is ${\cal B}[0,\infty)\otimes{\cal F}$-measurable in $(t,\omega')$ due to continuity of solution $\xi_t$ in $t$ and its measurability in $\omega'$ (see, e.g., \cite[Lemma 1.5.7]{Kry-ln}). Hence, the function $\widehat b(t,x,\omega'):=b(t,x,\xi_t(\omega'))$ is ${\cal B}[0,\infty)\otimes{\cal B}(\mathbb R^d)\otimes {\cal F}$- measurable in $(t,x,\omega')$. Further, one of the statements of Fubini theorem (see \cite[Theorem 1.5.5]{Kry-ln}) claims that
in this case the function
$$
\mathbb E' b(t,x,\xi_t) = \int b(t,x,\xi_t(\omega')) \mathbb P'(d\omega') = \int \widehat b(t,x,\omega') \mathbb P'(d\omega')
$$
is ${\cal B}[0,\infty)\otimes{\cal B}(\mathbb R^d)$-measurable, as required. Here we used the condition  of boundedness of $b$ in $y$ for each $x$, which implies integrability
$$
\iint\limits_{D} \int |b(t,x,\xi_t(\omega'))| \mathbb P'(d\omega') dtdx
\le \iint\limits_{D} C(x) dtdx
< \infty,
$$
over any {\em bounded} Borel subset $D \in {\cal B}[0,\infty)\otimes {\cal B}(\mathbb R^d).$     Lemma \ref{cor1} is proved.
\end{proof}

In this way, under the condition of the at most linear growth in $x$ assumed in the sequel (see (\ref{linear}) a few lines below), the coefficients $B$ and $\Sigma$ in (\ref{e1}) are Borel measurable in $(t,x)$; so, the equation (\ref{e1}) does make sense under this condition.

The next theorem is the main result of the paper about weak existence.

\begin{theorem}\label{thm1}
Let the initial value \(x_0\) have a finite 4th moment.
For the problem  (\ref{e1})-- (\ref{e200}), suppose that the following two conditions are  satisfied.
\begin{itemize}
\item[$(i)$] The functions $b$ and $\sigma$ admit
linear growth condition in $x$, i.e.,
there exists $C>0$ such that for any $s,x,y$,
\begin{equation}\label{linear}
|b(s,x,y)|+\|\sigma(s,x,y)\|\leq C(1+|x|),
\end{equation}

\item[$(ii)$] The diffusion matrix $\sigma$ is
uniformly non-degenerate in the following sense: there is a value $\nu>0$ such that for any probability measure $\mu$,
\begin{equation}\label{si}
\inf\limits_{s,x}\inf\limits_{|\lambda|=1}
\lambda^*\left(\int\sigma(s,x,y)\mu(dy)\right) \left(\int\sigma^*(s,x,y)\mu(dy)\right)
\lambda
\ge \nu.
\end{equation}
\end{itemize}
Then the equation (\ref{e1}) has a weak solution, that is,
a solution on some probability space with a standard $d_1$-dimensional
Wiener process with respect to some filtration $({\cal F}_t, \, t\ge 0)$.

\end{theorem}

\begin{remark}
If
$d_1=d$ and  $\sigma$ is {\em symmetric
and positive definite,} then the assumption (\ref{si}) can be replaced by an equivalent but much easier one, which is  frequently in use,
\begin{equation}\label{si1}
\inf\limits_{s,x,y}\inf\limits_{|\lambda|=1}
\lambda^*\sigma(s,x,y)
\lambda
\ge \nu >0.
\end{equation}
The  intuitive meaning of the condition  (\ref{si}) in the simplest 1D (that is, with $d_1=d=1$) situation is that the diffusion coefficient is non-degenerate and cannot change sign for any fixed $(s,x)$ and varying $y$.
It is plausible that any moment of order \(2+\epsilon\) for \(x_0\) suffices for all the  statements (except for  Theorem \ref{thm5} under the linear growth condition on the drift where an exponential moment will be required), but this goal is not pursued here.
Under the additional assumption of boundedness of \(b\) and \(\sigma\), the 4th moment of the initial value \(x_0\) is not necessary and can be further relaxed.
 \end{remark}

The structure of the proof of  Theorem \ref{thm1} is such that first  the case of symmetric non-degenerate $\sigma$ is tackled, i.e., with $d_1=d$ and
under the assumption (\ref{si1}).
A motivation for this approach is that under the relaxed assumption (\ref{si}) and under the symmetry of $\sigma$ it is easy to find a  smoothing of this matrix function which would keep the non-degeneracy of the diffusion coefficient.
Because of this we formulate some provisional statement under more restrictive  conditions as   Proposition \ref{pro22}. Its proof will be simultaneously the beginning of the proof of  Theorem~\ref{thm1}.

\medskip

\medskip

In the following intermediate simplified version of the Theorem \ref{thm1}, we allow both coefficients to grow, but $\sigma$ is assumed to be symmetric and positive definite. In the first part of the proof the coefficients are assumed to be bounded. We will treat now a more general structure of the coefficients $B$ and $\Sigma$, namely, we assume that instead of  $B$ and $\Sigma$ we have in \eqref{e1} coefficients of the form
\begin{equation*}
 \overline{ \Sigma}[t,X,\mu] = \phi\left(\int
\sigma(t,x,y)\mu(dy)\right), \quad
  \overline{B}[t,X,\mu] = \psi\left(\int
b(t,x,y)\mu(dy)\right),
\end{equation*}
where conditions on the matrix-valued functions $\phi: \mathbb R^{d\times d_1} \mapsto \mathbb R^{d\times d_1}$ (this includes the case of $d_1=d$ as in the Proposition \ref{pro22}) and vector-valued function $\psi: \mathbb R^{d} \mapsto \mathbb R^{d}$ will be specified. In fact, for the proof of Theorem \ref{thm1} we only need $\phi$; however, $\psi$ is added just by analogy since it does not bring any new difficulty. We use notations
$$\Sigma[t,X,\mu] = \int
\sigma(t,x,y)\mu(dy)=\langle \sigma, \mu\rangle_{t,x}, \;\;
B[t,X,\mu] = \int
b(t,x,y)\mu(dy)=\langle b, \mu\rangle_{t,x},
$$
and denote also
$$
A[t,x,\mu]:= \Sigma\Sigma^*[t,x,\mu].
$$
Actually, the conditions on the functions $\phi$ and $\psi$ imposed in the Proposition \ref{pro22} below might be further relaxed, but the authors do not have a motivation for that at the moment, as the main goal in this section is Theorem \ref{thm1}. Correspondingly, we now consider the equation \eqref{e1} with coefficients $  \overline{B}$ instead of $B$ and with $\overline{ \Sigma}$ instead of $\Sigma$.

\begin{Proposition}\label{pro22}
Assume that $d_1= d$,   let $\sigma$ and $b$ satisfy the linear growth condition (\ref{linear}), and for functions
$\phi$ and $\psi$ there exist $m_1, m_2>0$ such that for any $\Sigma', \Sigma'' \in \mathbb R^{d\times d}$ and any $B', B'' \in \mathbb R^{d}$
\begin{align}\label{fifi}
\|\phi(\Sigma') - \phi (\Sigma'')\|
\le C \|\Sigma' -  \Sigma''\|^{}(1+(\|\Sigma'\|\vee \|\Sigma''\|)^{m_1}),
 \\ \nonumber \\
\label{psipsi}
|\psi(B') - \psi(B'')|
\le C |B' - B''|(1+(|B'|\vee |B''|)^{m_2}),
\end{align}
as well as
\begin{align*}
\|\phi(\Sigma')\|
\le C (1+\|\Sigma'\|)^{},
 \\ \nonumber \\
|\psi(B')|
\le C (1+|B'|)^{},
\end{align*}
Also, $\overline{ \Sigma}[t,X,\mu]$ is assumed to be a symmetric positive definite uniformly non-degenerate for any $\mu$ and uniformly w.r.t.  $t,x,\mu$.

Let the initial value \(x_0\) have a finite fourth moment,  matrix $\sigma(t,x,y)$ be symmetric for each triple $(t,x,y)$, and the inequality (\ref{si1}) hold.
Then   equation (\ref{e1}) with coefficients $  \overline{B}$ and $\overline{ \Sigma}$ has a weak solution, that is,
a solution on some probability space with a standard $d$-dimensional
Wiener process with respect to some filtration $({\cal F}_t, \, t\ge 0)$.

\end{Proposition}

\begin{remark}
In the case of bounded coefficients, we use in the proof of Proposition \ref{pro22} the fact that the class of functions   $\phi$ and $\psi$ satisfying conditions \eqref {fifi} and \eqref{psipsi} must include identical $\psi=Id$
 and $\phi(\langle \sigma, \mu\rangle_{t,x}) = \sqrt{A[t,x,\mu]}$ defined via the Cauchy-Riesz-Dunford  formula for a function of a positive self-adjoint square root of the matrix $A$  (see, e.g., \cite[VII.3.9]{Dunford}),

\begin{equation}\label{contour}
\phi(\langle \sigma, \mu\rangle_{t,x}) = (2\pi i)^{-1}\oint_\Gamma \lambda^{1/2}(\lambda I- A[t,x,\mu])^{-1}d\lambda
\end{equation}
with $\Gamma = \Gamma(x) = \{\lambda: |\lambda| = r(x)\}$, whose radius $r(x)$ is such that this contour contains all eigenvalues of $A[t,x,\mu]$. For example, it can be $r(x) = \sup_{t,\mu}\|A[t,x,\mu]\|+1$. Under the boundedness condition on all coefficients which will be accepted temporary for the Proposition \ref{pro22} in the first part of the proof, it suffices to choose a unique contour with $\bar  r = \sup_{t,x,\mu}\|A[t,x,\mu]\|+1$, since this value is finite.

Now, let us check condition \eqref{fifi}
 in this  example with bounded coefficients. Indeed, we have for any probability measures $\mu, \mu'$ and any $t,x$,
\begin{align*}
\|\phi(\langle \sigma, \mu\rangle_{t,x})
- \phi(\langle \sigma, \mu'\rangle_{t,x})\|
  \nonumber\\ \nonumber\\
= \frac1{2\pi}\left\|\oint\limits_\Gamma \lambda^{1/2}(\lambda I- A[t,x,\mu])^{-1}d\lambda  - \oint\limits_\Gamma \lambda^{1/2}(\lambda I - A[t,x,\mu'])^{-1}d\lambda\right\|
  \nonumber\\ \nonumber\\
\le  C\bar r^{3/2} \sup_{|\lambda|=\bar r}\sup_{\nu} \|(\lambda I- A[t,x,\nu])^{-2}\|\,\|A[t,x,\mu] - A[t,x,\mu']\|
 \nonumber\\ \nonumber \\
\le C \bar r^{2} \,\|\Sigma[t,x,\mu] - \Sigma[t,x,\mu']\|,
\end{align*}

Recall that here $\phi(\langle \sigma, \mu\rangle_{t,x})= \phi(\mathbb E^3 \sigma(t,x,\xi))$ is assumed to be symmetric positive definite non-degenerate uniformly with respect to $t, x$ and $\mu$, where $\mu = {\cal L}(\xi)$.

Under the linear growth assumptions of Theorem \ref{thm1}, we will need $m_1=2$ and $m_2 = 0$, as we need to include
 $\psi = Id$ and   $ \phi(\langle \sigma, \mu\rangle)$ of the form \eqref{contour}
with $\Gamma = \Gamma(x) = \{\lambda: |\lambda| = r(x)\}$, whose  radius $r(x)$ is such that this contour contains all eigenvalues of $A[t,x,\mu]$.   It can be $r(x) = C|x|^2+1$ with some $C$ large enough.
In this case we have
\begin{align*}
\|\phi(\langle \sigma, \mu\rangle_{t,x})
- \phi(\langle \sigma, \mu'\rangle_{t,x})\|
  \nonumber\\ \nonumber\\
= \frac1{2\pi}\left\|\oint\limits_\Gamma \lambda^{1/2}(\lambda I- A[t,x,\mu])^{-1}d\lambda   - \oint\limits_\Gamma \lambda^{1/2}(\lambda I- A[t,x,\mu'])^{-1}d\lambda\right\|
  \nonumber\\ \nonumber\\
\le  C r(x)^{3/2} \sup_{|\lambda|=r(x)}\sup_{\nu} \|(\lambda I- A[t,x,\nu])^{-2}\|\,\|A[t,x,\mu] - A[t,x,\mu']\|
  \nonumber\\ \nonumber\\
=  C r(x)^{3/2} \sup_{|\lambda|=r(x)}\sup_{\nu} \|(\lambda I- A[t,x,\nu])^{-2}\|\,\|\Sigma\Sigma^*[t,x,\mu] - \Sigma\Sigma^*[t,x,\mu']\|
 \nonumber\\ \nonumber \\
\le C r(x)^{2} \|\Sigma[t,x,\mu] - \Sigma[t,x,\mu']\|.
\end{align*}
 Note that another slightly different option to define $\sqrt{A[t,x,\mu]}$ using the same idea in the unbounded case is as follows:
 \begin{equation}\label{sqrta2}
\sqrt{A[t,x,\mu]} = \frac1{2\pi i} \sum_{i=1}^{\infty}\mathbf{1}(i-1\le |x| < i) \oint_{\Gamma_i} \lambda^{1/2}(\lambda I- A[t,x,\mu])^{-1}\,d\lambda,
\end{equation}
where
$$
\Gamma_i = \{\lambda\in \mathbb C: \, |\lambda| = \sup_{t,x,\mu: \, |x|\le i}\|A[t,x,\mu]\| +1\},
$$
where the contour $\Gamma_i \subset \mathbb C$ in the complex plane is chosen in a way so that its interior contains  all the eigenvalues of the  elliptic  matrix $A[s,x,\cdot]$ for $|x|\le i$.

\end{remark}


\subsection{Proof of Proposition \ref{pro22}}
{\bf 1.} Assume temporarily that $\sigma$ and $b$ are bounded, and instead of (\ref{fifi})--(\ref{psipsi}) suppose for the first several steps that
\begin{equation}\label{phiA0}
\|\phi(\Sigma') - \phi(\Sigma'')\|\
\le C \|\Sigma' - \Sigma''\|^{},
\end{equation}
and
\begin{equation}\label{psiA0}
|\psi(B') - \psi(B'')|
\le C |B'' - B''|^{}.
\end{equation}
These restrictions will be waived at the step 6 of the proof.
Let us  smooth out both coefficients with respect to all variables by convolutions in such a way that they become globally Lipschitz in $x$ and $y$. Namely, let
\begin{equation*}
b^n(t,x,y) = b(t,x,y) * \varphi_n(x) *
\varphi_n(y),
\end{equation*}
and
\begin{equation*}
\sigma^n(t,x,y) = \sigma(t,x,y) * \varphi_n(x) *
\varphi_n(y),
\end{equation*}
 where the sequence $\varphi_n(\cdot)$ is defined in a
standard way, i.e., as non-negative \(C^\infty\) functions with a compact support integrated to one, and so that this compact support squeezes to the origin of the corresponding variable as \(n\to\infty\); or, in other words, that they are delta-sequences in the corresponding variables.
Note that, of course, we may assume that for every $n$ the smoothed coefficient of the drift remains to be under the linear growth condition (\ref{linear}) with the same constant for each $n$ (in reality this constant may increase a little bit in comparison to the constant $C$ from (\ref{linear}), but still remain uniformly bounded); also, under the assumption (\ref{si1}) the smoothed diffusion remains uniformly non-degenerate with ellipticity constants independent of~$n$.

\medskip

Now we shall explain why the equation with smoothed coefficients has a (strong) solution.
We use successive approximations. For any fixed $n$, let
$$
X^n(0)_t:= x_0, \quad \mu^n(0)_t={\cal L}(X^n(0)_t) = \mu_0, \quad t\ge 0;
$$
further, if $X(m)_t$ and $\mu(m)_t$ are already determined, let us define
$$
X^n(m+1)_t:= x_0 + \int_0^t \bar B^n[s,X^n(m)_s, \mu^n(m)_s]ds + \int_0^t  \bar \Sigma^n[s,X^n(m)_s, \mu^n(m)_s]dW_s,
$$
where
$$
\bar B^n[t,x,\mu] = \psi\left(\int
b^n(t,x,y)\mu(dy)\right), \;\;
\bar \Sigma^n[t,x,\mu] =\phi\left(\int
\sigma^n(t,x,y)\mu(dy)\right).
$$
We can say that $\bar  B^n[s,X^n(m)_s, \mu^n(m)_s] = \psi(\mathbb E^3 b^n(s,X^n(m)_s, \xi^n(m)_s))$, where $\xi^n(m)_s$ is a random variable equivalent to $X^n(m)_s$ on some independent probability space, and, moreover, the sequence $(\xi^n(m)_\cdot), m\ge 1$ can be chosen independent on $(X^n(m)_\cdot), m\ge 1$, and so that the whole sequence $(\xi^n(m)_\cdot), m\ge 1$ has the same distribution as the sequence $(X^n(m)_\cdot), m\ge 1$.
Then by induction the second moments of any $X^n(m)_t$ are finite and uniformly bounded for $t\le T$, and by It\^o's isometry and CBS inequality,
\begin{align*}
\mathbb E |X^n(m+1)_t - X^n(m)_t|^2
\\
\le C_T \mathbb E \int_0^t |\psi(\mathbb E^3 b^n(s,X^n(m)_s, \xi^n(m)_s))
- \psi(\mathbb E^3 b^n(s,X^n(m-1)_s, \xi^n(m-1)_s))|^2ds
\\
+ C  \mathbb E \int_0^t \|\phi(\mathbb E^3 \sigma^n(s,X^n(m)_s, \xi^n(m)_s))
- \phi(\mathbb E^3\sigma^n(s,X^n(m-1)_s, \xi^n(m-1)_s))\|^2ds
\\
\le C_T \mathbb E \int_0^t |\mathbb E^3 b^n(s,X^n(m)_s, \xi^n(m)_s)
- \mathbb E^3 b^n(s,X^n(m-1)_s, \xi^n(m-1)_s)|^2ds
\\
+ C  \mathbb E \int_0^t \|\mathbb E^3 \sigma^n(s,X^n(m)_s, \xi^n(m)_s)
- \mathbb E^3\sigma^n(s,X^n(m-1)_s, \xi^n(m-1)_s)\|^2ds
\\
\le C_{T,n} \mathbb E \int_0^t |X^n(m)_s - X^n(m-1)_s|^2ds
+ C  \mathbb E \int_0^t \mathbb E^3 |\xi^n(m)_s
- \xi^n(m-1)_s|^2ds
\\
\le C_{T,n} \mathbb E \int_0^t |X^n(m)_s - X^n(m-1)_s|^2  ds.
\end{align*}
Since all terms here are finite, we obtain by induction
\begin{align*}
\mathbb E |X^n(m+1)_t - X^n(m)_t|^2
\le  C_{T,n}  \frac{T^m}{m!}, \quad t\le T.
\end{align*}
Due to the Doob inequality we also get
\begin{align*}
\mathbb E \sup_{t\le T}|X^n(m+1)_t - X^n(m)_t|^2 \le C_{T,n}  \frac{T^m}{m!}.
\end{align*}
From here  by standard methods it follows easily convergence  in probability of the sequence $X^n(m)_t$ to a solution $X^n_t$ as $m\to\infty$, uniformly with respect to $t\le T$, as required.

\medskip

\noindent
{\bf 2}. In a standard way (see, e.g.,  \cite{Kry}, \cite{Sko}), absolutely similar to the inequalities of Lemma \ref{lebds},  we get the estimates uniform in $n$,
\begin{align}
\label{kol1}
 \mathbb{E} \sup_{0\le t\le T} |X_t^n |^4 \le C_T (1+\mathbb E |x_0|^4),
\end{align}
and
\begin{eqnarray}\label{kol2}
\sup_{0\le s\le t\le T; \, t-s\le h}\mathbb{E} |X_t^n - X^n_s|^4 \le
C_{T}  h^2,
\end{eqnarray}
with some constants $C_{T}$ which may be different for different inequalities but do not depend on $n$.
In fact, similar a priori bounds hold true for any power function assuming the appropriate initial moment, although, this will not be used in this paper. The proof can be done following the lines in  \cite[Theorem 1.6.4]{GS68}.
Why do we need the fourth degree will be clear in the next step: it is useful for verifying continuity for the processes with equivalent finite-dimensional distributions. Note that all these bounds are valid under the linear growth assumptions in $x$.

\medskip

\noindent
{\bf 3}.
 Let us introduce new processes $\xi^n$, the copies of  $X^n$, on some other  independent   probability space (i.e., we will consider both on the direct product of the two probability spaces); it also satisfies a similar SDE. In the sequel by $ \mathbb{E}^3  \sigma^n(s,X^n_s,\xi^n_s) $ or  $ \mathbb{E}^3  \sigma(s,X_s,\xi_s) $ we denote
expectation with respect to the third variable $ \xi^n_s $, or $ \xi_s $ i.e.,
{\em conditional} expectation given the second variable $ X^n_s$ or $X_s$; in other words, $$ \mathbb{E}^3  \sigma^n(s,X^n_s,\xi^n_s) = \int \sigma^n(s,X^n_s,y)\mu^{\xi^n}_s(dy),$$  where $ \mu^{\xi^n}_s $ stands for the marginal distribution of $ \xi^n_s $; likewise,
$$
\mathbb{E}^3  (\sigma^n(s,X^n_s,\xi^n_s)
- \sigma^n(s,X_s,\xi_s))
$$
means simply
$$
\displaystyle \int \sigma^n(s,X^n_s,y)\mu^{\xi^n}_s(dy)
-  \int \sigma^n(s,X_s,y)\mu^{\xi}_s(dy),
$$
where  $\mu^{\xi}_s$ is the marginal distribution of $\xi_s$,
and, finally,
$$
\mathbb{E}^3 \|\sigma^n(s, X^n_s, \xi^n_s) - \sigma(s,X_s, \xi_s)\|^2
$$
is understood as
$$
\displaystyle \int \|\sigma^n(s,X^n_s,y)
-  \sigma^n(s,X_s,y')\|^2\mu^{\xi^n,\,\xi}_s(dy,dy'),
$$
where $\mu^{\xi^n,\,\xi}_s(dy,dy')$ denotes the marginal distribution of the couple $(\xi^n_s,\xi_s)$.

\medskip

Now, due to the estimates  (\ref{kol1})--(\ref{kol2}) and by virtue of Skorokhod's
Lemma about a single probability
space and convergence in probability (see Lemma \ref{app2} in the Appendix, or \cite[\S 6, ch. 1]{Sko},
or \cite[Lemma 2.6.2]{Kry},
without loss of generality we  assume that
not only $\mu^n \Longrightarrow \mu$, but also on some probability space
$$
(\widetilde X^n_{t},\widetilde \xi^n_{t}, \widetilde W^n_{t}) \stackrel{\mathbb{P}}{\to}
(\widetilde X_t,\widetilde \xi_t, \widetilde W_t), \quad n\to\infty,
$$
 for any~$t$ and  for some   equivalent  random processes $(\widetilde X^n, \widetilde \xi^n, \widetilde W^n)$,
generally speaking, over a sub-sequence. Slightly abusing notations, we  denote initial values still by \(x_0\) without tilde.
Also,  without loss of generality we assume that each process $(\widetilde \xi ^n_{t}, \, t\ge
0)$ for any $n\ge 1$ is independent of $(\widetilde X^n, \widetilde W^n)$, as well as their limit $\widetilde \xi_t$ may be chosen  independent
of the limits $(\widetilde X,\widetilde W)$ (this follows from the fact that on the original probability space $\xi^n$ is independent of $(X^n, W^n)$ and on the new probability space their joint distribution remains the same; hence, independence of $\widetilde \xi^n$ is also valid and in the limit this is still true). See the details in the proof of the Theorem 2.6.1 in \cite{Kry}. We may
also introduce Wiener processes for $\xi^n_{t}$ and $\widetilde \xi^n_{t}$, and will do it because it will be useful at one of the steps of the proof of the Proposition \ref{thm4}  (for the Proposition \ref{pro22} it is not necessary). Namely, on an independent probability spaces we have,
\begin{equation}\label{exi0}
d\xi^n_t = \overline B^{n}[t,\xi^n_t, \mu_t]dt +  \overline \Sigma^{n} [t,\xi^n_t, \mu_t]dW^{\prime,n}_t, \;\; t\ge 0, \qquad
{\cal L}(\xi^n_0)={\cal L}(x_0),
\end{equation}
and
\begin{equation*}
d\widetilde \xi^n_t =  \overline B^{n}[t,\widetilde \xi^n_t, \mu_t]dt +  \overline \Sigma^{n} [t,\widetilde \xi^n_t, \mu_t]d\widetilde W^{\prime,n}_t, \;\; t\ge 0, \qquad
{\cal L}(\widetilde \xi^n_0)={\cal L}(x_0).
\end{equation*}

In what follows, let us fix some arbitrary \(T>0\) and consider \(t\) in the interval \([0,T]\).
Due to the inequality (\ref{kol2}),
the same inequality holds for $\widetilde X^n$ and  $\widetilde W^n$, in particular,
\begin{eqnarray}\label{kol2t}
\sup_{0\le s\le t\le T; \, t-s\le h}\mathbb{E} |\widetilde X_t^n - \widetilde X^n_s|^4 \le
C_{T}  h^2.
\end{eqnarray}
Due to Kolmogorov's continuity theorem, it means that all processes $\widetilde X^n$ may be regarded as continuous, and  $\widetilde W^n$ can be assumed also continuous by the same reason.
Further, due to the independence of the increments of $W^n$ after time $t$ of the sigma-algebra $\sigma(X^n_s, W^n_s, s\le t)$, the same property holds true for $\widetilde W^n$ and  $\sigma(\widetilde X^n_s, \widetilde W^n_s, s\le t)$, as well as for $\widetilde W^n$ and for the completions of the sigma-algebras  $\sigma(\widetilde X^n_s, \widetilde W^n_s, s\le t)$ which we denote by ${\cal F}^{(n)}_t$. Also,  the processes  $\widetilde X^n$ are adapted to the filtration $({\cal F}^{(n)}_t)$. So, all stochastic integrals which involve
$\widetilde X^n$ and $\widetilde  W^n$ are well defined. The same relates to the processes $\widetilde \xi^n$.

Hence,
again by using Skorokhod's Lemma (Lemma \ref{app2}),  we may choose a
subsequence $n'\to\infty$ so as to pass to the limit in the equation
$$
\widetilde X^{n'}_t = x_0 + \int_0^t \psi(\mathbb{E}^3 b^{n'}(s,\widetilde X^{n'}_s, \widetilde \xi^{n'}_s))\,ds +
\int_0^t \phi(\mathbb{E}^3 \sigma^{n'}(s,\widetilde X^{n'}_s,  \widetilde \xi^{n'}_s)) d\widetilde W^{n'}_s,
$$
in order to get
$$
\widetilde X_t = x_0 + \int_0^t \psi(\mathbb{E}^3 b(s,\widetilde X_s, \widetilde \xi_s))ds +
\int_0^t  \phi(\mathbb{E}^3 \sigma(s,\widetilde X_s, \widetilde  \xi_s))\, d\widetilde W_s,
$$
or, equivalently,
$$
\widetilde X_t = x_0 + \int_0^t  \bar {B}[s,\widetilde X_s, \mu_s]\,ds + \int_0^t   \bar {\Sigma}[s,\widetilde X_s,  \mu_s]\, d\widetilde W_s,
$$
with
$$
\mu_s = {\cal L}(\widetilde X_s).
$$

First of all, recall that a priori bounds
(\ref{kol1}) -- (\ref{kol2}) and (\ref{kol2t}) hold true with constants not depending on $n$.
Now, by Skorokhod's Lemma   \ref{app2}, on some probability space we have a sequence of  equivalent processes $(\widetilde X^{n'}_t,\widetilde \xi^{n'}_t, \widetilde W^{n'}_t)$ and a limiting triple $(\widetilde X^{}_t,\widetilde \xi^{}_t, \widetilde W^{}_t)$ such that for any $t$,
\begin{equation*}
(\widetilde X^{n'}_t,\widetilde \xi^{n'}_t, \widetilde W^{n'}_t) \stackrel{\mathbb P}{\to} (\widetilde X^{}_t,\widetilde  \xi^{}_t, \widetilde W^{}_t).
\end{equation*}
By virtue of the above a priori estimates for \(\widetilde W^n\), the process \(\widetilde W\) is continuous and it is a Wiener process. 

Also, using the properties of the triple $(\widetilde X^{n'}_t,\widetilde \xi^{n'}_t, \widetilde W^{n'}_t)$ before the limit, it may be verified that the process $\widetilde W_t$ is a Wiener one with respect to the filtration $\widetilde {\cal F}_t:= {\cal F}^{\tilde X, \tilde W}_t$. Namely, the increments of  $\widetilde W_t$, say, after $t_0$ are independent on the sigma-algebra ${\cal F}^{\tilde X, \tilde W}_{t_0}$, or, in other words, sigma-algebras $\sigma(\widetilde W_t-\widetilde W_s: \, t_0\le s\le t)$ and ${\cal F}^{\tilde X, \tilde W}_{t_0}$ are independent for any $t_0$. This follows straightforwardly from the same assertion before the limit.

In particular,  related   Lebesgue and stochastic integrals are all well defined.
Moreover, by virtue of the uniform estimates (\ref{kol2}), the limit \((\widetilde X^{}_t,\widetilde  \xi^{}_t)\) may be also regarded as  continuous due to Kolmogorov's continuity theorem  because   a priori bounds  (\ref{kol1}) -- (\ref{kol2}) remain valid for the limiting processes \(\widetilde X, \widetilde \xi\). In particular,   note because it will come in handy later,  that
\begin{eqnarray*}
\sup_{0\le t\le T} \mathbb{E} |\widetilde X_t|^2 \le C_T (1+\mathbb E |x_0|^2).
\end{eqnarray*}

\noindent
{\bf 4}.
\noindent
Still in the case of bounded coefficients $b$ and $\sigma$, let us now show that
\begin{equation}\label{limNb0}
\int_0^t \psi(\mathbb{E}^3 b^{n'}(s,\widetilde X^{n'}_s, \widetilde \xi^{n'}_s))ds  \stackrel{\mathbb P}{\to} \int_0^t \psi( \mathbb{E}^3 b(s,\widetilde X_s, \widetilde \xi_s))ds,
\end{equation}
and
\begin{equation}\label{limNs0}
\int_0^t \phi(\mathbb{E}^3 \sigma^{n'}(s,\widetilde X^{n'}_s,  \widetilde \xi^{n'}_s)) d\widetilde W^{n'}_s  \stackrel{\mathbb P}{\to} \int_0^t  \phi(\mathbb{E}^3 \sigma(s,\widetilde X_s, \widetilde  \xi_s)) d\widetilde W_s, \quad n'\to\infty.
\end{equation}
We start with the drift term.
Let us fix some  $n_0$ and let $n>n_0$. Due to the assumption (\ref{psiA0}), we have for any $t\le T$,
\begin{align*}
& \displaystyle \mathbb P \left(\left|\int_0^t \left(\psi(\mathbb{E}^3  b^{n}(s,\widetilde X^{n}_s, \widetilde \xi^{n}_s))  - \psi(\mathbb{E}^3  b(s,\widetilde X_s,\widetilde  \xi_s))\right)ds\right| > c\right)
 \nonumber \\\nonumber \\\nonumber
& \displaystyle \le \mathbb P\left(
\left|\int_0^t  \left(\psi(\mathbb{E}^3 b^{n}(s,\widetilde X^{n}_s, \widetilde \xi^{n}_s) ) - \psi(\mathbb{E}^3 b^{n_0}(s,\widetilde X^n_s,\widetilde  \xi^n_s))\right)ds\right| > \frac{c}{3}\right)
 \nonumber \\\nonumber \\\nonumber
& \displaystyle + \mathbb P \left(\left|\int_0^t  \left(\psi( \mathbb{E}^3 b^{n_0}(s,\widetilde X^{n}_s, \widetilde \xi^{n}_s))  - \psi( \mathbb{E}^3b^{n_0}(s,\widetilde X_s,\widetilde  \xi_s))\right)ds\right| > \frac{c}3\right)
 \nonumber \\\nonumber \\\nonumber
& \displaystyle + \mathbb P \left(\left|\int_0^t  \left(\psi( \mathbb{E}^3 b^{n_0}(s,\widetilde X^{}_s, \widetilde \xi^{}_s))  - \psi( \mathbb{E}^3 b(s,\widetilde X_s,\widetilde  \xi_s))\right)ds\right| > \frac{c}3\right)
  \\ \nonumber\\
& \displaystyle =: I^1 + I^2 + I^3.
\end{align*}
Let
\begin{align*}
\gamma_{n,R}:= \inf(t\ge 0: \sup_{0\le s\le t} (|\widetilde X_s^n|\vee  |\widetilde \xi_s^n|) \ge R),
 \;
\gamma_{R}:= \inf(t\ge 0: \sup_{0\le s\le t} (|\widetilde X_s|\vee  |\widetilde \xi_s|) \ge R),
 \\
\gamma^X_{R}:= \inf(t\ge 0: \, \sup_{0\le s\le t} |\widetilde X_s| \ge R), \;
\gamma^\xi_{R}:= \inf(t\ge 0: \, \sup_{0\le s\le t} |\widetilde \xi_s| \ge R),
 \\
\gamma^X_{n,R}:= \inf(t\ge 0: \, \sup_{0\le s\le t} |\widetilde X_s^n|\ge R),
\;
\gamma^\xi_{n,R}:= \inf(t\ge 0: \, \sup_{0\le s\le t} |\widetilde \xi_s^n| \ge R).
\end{align*}
We have that
for any $\epsilon>0$ there exists $R>0$ such that (to have \(R-1\) instead of \(R\) will be convenient shortly)
\[
\mathbb P(\sup_{0\le t\le T} (|\widetilde X_t|\vee  |\widetilde \xi_t|) \ge   R-1) < \epsilon,
\]
or, equivalently,
\begin{equation*}
\mathbb P(\gamma_{R-1}\le T)  < \epsilon,
\end{equation*}
and similarly,
\[
\sup_n \mathbb P(\sup_{0\le t\le T} (|\widetilde X^n_t|\vee  |\widetilde \xi^n_t|) \ge   R-1) < \epsilon,
\]
or, equivalently,
\begin{equation*}
\sup_n \mathbb P(\gamma_{n,R-1}\le T)  < \epsilon.
\end{equation*}


Denote
\[
g^{n,n_0}(s,x,\xi):=b^{n}(s,x,\xi)- b^{n_0}(s,x,\xi), \quad
g^{n_0}(s,x,\xi):=b^{n_0}(s,x,\xi)- b^{}(s,x,\xi).
\]
Then the first summand $I^1$  may be estimated by BCM   inequality   due to the condition (\ref{psiA0}) as follows without $\psi$,
\begin{eqnarray*}\label{i1}
I^1 \le \frac3{c} \, \mathbb E  \int_0^T C \,\mathbb{E}^3 |b^{n}(s,\widetilde X^{n}_s,\widetilde  \xi^{n}_s)- b^{n_0}(s,\widetilde X^n_s, \widetilde \xi^n_s)|\,ds
 \nonumber\\ \nonumber
=  \frac{3C}{c}\, \mathbb E \mathbb{E}^3 \int_0^T  |g^{n,n_0}(s,\widetilde X^{n}_{s},\widetilde  \xi^{n}_{s})|\,ds
=  \frac{3C}{c}\, \mathbb E \int_0^T  |g^{n,n_0}(s,\widetilde  X^{n}_{s},\widetilde   \xi^{n}_{s})|\,ds
\\
= C \mathbb E \mathbf{1}(\gamma_{n,R}\le T) \int_0^T  |g^{n,n_0}(s,\widetilde  X^{n}_{s},\widetilde   \xi^{n}_{s})|\,ds
+C \mathbb E\mathbf{ 1}(\gamma_{n,R} > T) \int_0^T  |g^{n,n_0}(s,\widetilde  X^{n}_{s},\widetilde   \xi^{n}_{s})|\,ds.
\end{eqnarray*}
Here the first term $\displaystyle \mathbb E \mathbf{1}(\gamma_{n,R}\le T) \int_0^T  |g^{n,n_0}(s,\widetilde  X^{n}_{s},\widetilde   \xi^{n}_{s})|\,ds$ admits the bound (recall that on each line constants $C$ can be different)
\begin{align*}
\mathbb E \mathbf{1}(\gamma_{n,R}\le T) \int_0^T  |g^{n,n_0}(s,\widetilde  X^{n}_{s},\widetilde   \xi^{n}_{s})|\,ds
\le CT^{} \mathbb P(\gamma_{n,R}\le T)^{},
\end{align*}
and the  last expression is small uniformly in $n$ if $R$ is large enough.

The second term $\mathbb E \mathbf{1}(\gamma_{n,R} > T) \int_0^T  |g^{n,n_0}(s,\widetilde  X^{n}_{s},\widetilde   \xi^{n}_{s})|\,ds$ admits a bound via   Krylov's estimate (see  the Theorems 2.4.1 or 2.3.4 in \cite{Kry}) as follows: there exists a constant $N$ depending on the dimension $d$ and on $R$ through the ellipticity constants of the diffusion matrix and the sup-norm of the drift on the set $B_R \times B_R$, where $B_R=\{x\in \mathbb{R}^d: |x|\le R\}$, so that
\begin{eqnarray*}
& \displaystyle\mathbb E \mathbf{1}(\gamma_{n,R} > T)\int_0^T  |g^{n,n_0}(s,\widetilde  X^{n}_{s},\widetilde   \xi^{n}_{s})|\,ds
\\
& \displaystyle\le N_R \left(\int_0^T \int_{|x|\le \widetilde R}\int_{|\xi|\le \widetilde R} |g^{n,n_0}(s,x,\xi)|^{2d+1}\,dx d\xi ds\right)^{\frac1{2d+1}}
\\
& \displaystyle \le
N_R\, \left(\int_0^T \int_{|x|\le \widetilde R}\int_{|\xi|\le \widetilde R} |b^{n}(s,x,\xi) - b^{}(s,x,\xi)|^{2d+1}\,dx d\xi ds\right)^{\frac1{2d+1}}
 \\
& \displaystyle + N_R\, \left(\int_0^T \int_{|x|\le \widetilde R}\int_{|\xi|\le \widetilde R} |b^{n_0}(s,x,\xi) - b^{}(s,x,\xi)|^{2d+1}\,dx d\xi ds\right)^{\frac1{2d+1}}
\to 0, \quad n, n_0 \to \infty,
\end{eqnarray*}
for each $R$, by virtue of the well-known property of mollified functions. Hence, overall, we obtain that
$$
I^1 \to 0, \quad n, n_0\to\infty.
$$

\medskip

Further, under the assumptions of   Proposition \ref{pro22}, the second term $I^2$ admits for any \(0\le t\le T\) the estimate, again without $\psi$ in the right hand side,
\begin{eqnarray*}
&\displaystyle \!\!\! I^2 \!=\!
\mathbb P \!\left(\left|\!\int_0^t \!\left(\psi(\mathbb{E}^3 b^{n_0}(s,\widetilde X^{n}_s, \widetilde \xi^{n}_s))  \! -\! \psi(\mathbb{E}^3 b^{n_0}(s,\widetilde X_s,\widetilde  \xi_s))\right)\!ds\right| \!>\! \frac{c}3\right)
\\
&\displaystyle \le C\, \mathbb E \mathbb{E}^3 \int_0^T  |b^{n_0}(s,\widetilde X^{n}_s, \widetilde \xi^{n}_s)- b^{n_0}(s,\widetilde X_s, \widetilde \xi_s)|\,ds,
\end{eqnarray*}
which
tends to zero as $n\to \infty$ under the assumptions of the Proposition \ref{pro22}, due to the Lebesgue bounded convergence theorem because of the convergence in probability $(\widetilde X^{n}_s, \widetilde \xi^{n}_s) \to (\widetilde X^{}_s, \widetilde \xi^{}_s)$. So, for each $n_0$
$$
\lim_{n\to\infty} I^2 = 0.
$$

\medskip

As for the   term $I^3$, it admits the following bounds (again without $\psi$ in the right hand side):
\begin{eqnarray*}
\!\!& \displaystyle I^3\!
=\! \mathbb P \!\left(\left|\!\int_0^t \!\left(\psi(\mathbb{E}^3 b^{n_0}(s,\widetilde X^{}_s, \widetilde \xi^{}_s)) \! - \! \psi(\mathbb{E}^3 b^{}(s,\widetilde X_s,\widetilde  \xi_s))\right)\!ds\right|\! > \! \frac{c}3\right)
 \nonumber \\ \nonumber \\ \nonumber
& \displaystyle \le C\, \mathbb E \mathbb{E}^3\int_0^T  |b^{n_0}(s,\widetilde X^{}_s, \widetilde \xi^{}_s)- b^{}(s,\widetilde X_s, \widetilde \xi_s)|\,ds
 \nonumber \\ \nonumber \\
& \displaystyle =C\,\mathbb E
(\mathbf{1}(\gamma_{n,R}\wedge \gamma_{R} \le T) + \mathbf{1}(\gamma_{n,R}\wedge \gamma_{R} > T))
\int_0^T |g^{n_0}(s,\widetilde X^{}_s, \widetilde \xi^{}_s)|\,ds.
\end{eqnarray*}

Furthermore, the term $\displaystyle \mathbb E
\mathbf{1}(\gamma_{n,R}\wedge \gamma_{R} \le T)
\int_0^T |g^{n_0}(s,\widetilde X^{}_s, \widetilde \xi^{}_s)|\,ds$ admits the bound
\begin{align*}
\mathbb E
\mathbf{1}(\gamma_{n,R}\wedge \gamma_{R} \le T)
\int_0^T |g^{n_0}(s,\widetilde X^{}_s, \widetilde \xi^{}_s)|\,ds
\le C \mathbb P (\gamma_{n,R}\wedge \gamma_{R} \le T) \to 0, \quad R\to \infty,
\end{align*}
uniformly in $n$.
Hence, it follows that
$$
\lim_{n_0 \to\infty} \lim_{R \to\infty}\mathbb E
\mathbf{1}(\gamma_{n,R}\wedge \gamma_{R} \le T)
\int_0^T |g^{n_0}(s,\widetilde X^{}_s, \widetilde \xi^{}_s)|\,ds = 0.
$$

In order to evaluate the term $\mathbb E
\mathbf{1}(\gamma_{n,R}\wedge \gamma_{R} > T)
\int_0^T |g^{n_0}(s,\widetilde X^{}_s, \widetilde \xi^{}_s)|\,ds$, we note that
the values of the function $g^{n_0}(s,x,\xi)$ outside the set $\{(x,\xi):\, (|x|\vee |\xi|)\le   R\}$ with large $  R$
whose exact value is  not relevant. So, for evaluating this term, without losing of generality we may assume  that $g^{n_0}(s,x,\xi)$ vanishes outside the ball $B_{R+1}\times B_{R+1}$: if not, we just truncate accepting that $g^{n_0}=0$ outside $B_{R+1}\times B_{R+1}$. Then the desired convergence follows from  Krylov's bound of Lemma  \ref{lekrybd}. Indeed,
provided that we denote $$g_R^{n_0}(s,x,\xi):= g^{n_0}(s,x,\xi)\mathbf{ 1}(|x|\le R, |\xi|\le R),$$
the following bounds follow immediately:
\begin{align*}
\mathbb E
\mathbf{1}(\gamma_{n,R}\wedge \gamma_{R} > T)
\int_0^T |g^{n_0}(s,\widetilde X^{}_s, \widetilde \xi^{}_s)|\,ds
\le \mathbb E
\int_0^T |g_R^{n_0}(s,\widetilde X^{}_s, \widetilde \xi^{}_s)|\,ds.
\end{align*}
Now,
for continuous $g^{n_0}$, in the limit as $n\to\infty$  we get
\begin{align*}
\mathbb E
\mathbf{1}( \gamma_{R} > T)
\int_0^T |g^{n_0}(s,\widetilde X^{}_s, \widetilde \xi^{}_s)|\,ds
\le \mathbb E
\int_0^T |g_R^{n_0}(s,\widetilde X^{}_s, \widetilde \xi^{}_s)|\,ds
\\
= \lim_{n\to\infty}\mathbb E
\int_0^T |g_R^{n_0}(s,\widetilde X^{n}_s, \widetilde \xi^{n}_s)|\,ds
\\
\le C N_R\|g_R^{n_0}\|_{L_{2d+1}((0,\infty)\times \mathbb R^d\times \mathbb R^d)}
= C N_R\|g_R^{n_0}\|_{L_{2d+1}((0,\infty)\times B_R\times B_R)}.
\end{align*}
The latter bound extends 
to all Borel measurable functions  $g_R$ as in  \cite[Section II]{Kry}.

\begin{equation}\label{kb}
\mathbb E
\int_0^T  |g_R^{}(s,\widetilde X^{}_s, \widetilde  \xi^{}_s)|\,ds \le N_R \|g_R\|_{L_{2d+1}([0,T]\times B_R \times B_R)}.
\end{equation}

Let us show that this bound remains true for any Borel function $g_R$ vanishing outside $B_R \times B_R$. This is analogous to the justification of a similar inequality in \cite[Section II.6]{Kry}, but we want to add some details for the convenience of the reader. Firstly,  it clearly suffices to prove it for non-negative functions only, $g_R\ge 0$. Secondly, since indicator of any compact in a finite-dimensional Euclidean space may be approximated monotonically (from above) pointwise by continuous bounded functions, again by Fatou lemma  (\ref{kb}) remains valid for all non-negaitive functions represented by finite sums
$$
g(s,x,\xi) = \sum_{i=1}^{K} c_i 1_{\Gamma_i}(s,x,\xi),
$$
with any compacts $\Gamma_i \subset [0,T]\times B_R \times B_R$ and with all $c_i>0$ by the monotonic convergence theorem. Thirdly, the expression $\mathbb E
\int_0^T  g^{}(s,\widetilde X^{}_s, \widetilde  \xi^{}_s)\,ds$
 may be rewritten via a ``Green measure''
$$
\mathbb E
\int_0^T  g^{}(s,\widetilde X^{}_s, \widetilde  \xi^{}_s)\,ds
=
\int g^{}(s,x, \xi)\,\nu(ds, dx, d\xi)
$$
with a positive sigma-additive measure $\nu$:
$$
\nu(ds, dx, d\xi) = ds\, \P(\tilde X_s \in dx)\P(\tilde \xi_s \in d\xi).
$$ Fourthly, since on finite-dimensional Euclidean spaces any such measure is regular, it is possible to approximate monotonically any expression of the form
$$
\int g^{}(s,x, \xi)\,\nu(ds, dx, d\xi) =
\int (\sum_{i=1}^{K} c_i 1_{D_i}(s,x,\xi))\,\nu(ds, dx, d\xi)
$$
with any Borel $D_i \subset [0,T]\times B_R \times B_R$ and $c_i>0$  by  $\int (\sum_{i=1}^{K} c_i 1_{\Gamma_i}(s,x,\xi))\,\nu(ds, dx, d\xi)$ with compact $\Gamma_i$. Hence, for such functions $\sum_{i=1}^{K} c_i 1_{D_i}(s,x,\xi)$ with Borel $D_i$ the estimate (\ref{kb}) also remains valid.
Finally, by such sums $\sum_{i=1}^{K} c_i 1_{D_i}(s,x,\xi)$ any bounded Borel measurable $g_R\ge 0$ can be approximated uniformly and monotonically; so, the bound  (\ref{kb}) is valid for any bounded Borel measurable $g_R\ge 0$, as required. For what follows, let us note that it also suffices for extension of the estimate to all Borel measurable from $L_{2d+1}([0,T]\times B_R \times B_R)$ again by monotonic approximations.

Hence, by the properties of the mollified functions  it follows that
$$
\lim_{n_0\to\infty}I^3 = 0.
$$
The convergence  (\ref{limNb0}) is, thus, proved.

\medskip

\noindent
{\bf 5}. Now still for bounded coefficients let us consider convergence of stochastic integrals in (\ref{limNs0}). Our goal is an estimate similar to that for the drift and Lebesgue integrals above:
\begin{equation}\label{sest0}
\mathbb P\left(\left\|\int_0^t \phi(\mathbb{E}^3 \sigma^n(s,\widetilde X^n_s,  \widetilde \xi^n_s)) d\widetilde W^n_s  - \int_0^t  \phi(\mathbb{E}^3 \sigma(s,\widetilde X_s,  \widetilde \xi_s)) d\widetilde W_s\right\| > c\right)
 < C \epsilon,
\end{equation}
for any $c, \epsilon>0$ and $n$ large enough. In principle, the task is similar to the convergence of Lebesgue integrals studied in the previous steps. Hence, we mainly show how to overcome the additional obstacle due to different Wiener processes \(\widetilde W \) and \(\widetilde W^n\) in the stochastic integrals. We have a tool for this which is Skorokhod's Lemma  \ref{app1}.

By virtue of \cite[Theorem 6.2.1(v)]{Kry-ln} and similarly to the calculus for the drift with Lebesgue integrals in the previous steps, yet using second moments instead of the first ones for the evident reason, we estimate
\begin{align*}
& \displaystyle \mathbb P\left(\left\|\int_0^t \phi(\mathbb{E}^3 \sigma^n(s,\widetilde X^n_s,  \widetilde \xi^n_s)) d\widetilde W^n_s  - \int_0^t  \phi(\mathbb{E}^3 \sigma(s,\widetilde X_s,  \widetilde \xi_s)) d\widetilde W_s\right\| > c\right)
 \\
& \displaystyle \le \mathbb P\left(\left\|\int_0^t \phi(\mathbb{E}^3 \sigma^n(s,\widetilde X^n_s,  \widetilde \xi^n_s)) d\widetilde W^n_s  - \int_0^t  \phi(\mathbb{E}^3 \sigma ^{n_0}(s,\widetilde X^n_s,  \widetilde \xi^n_s)) d\widetilde W^n_s\right\| > c/3\right)
 \\
& \displaystyle + \mathbb P\left(\left\|\int_0^t \phi(\mathbb{E}^3 \sigma^{n_0}(s,\widetilde X^n_s,  \widetilde \xi^n_s)) d\widetilde W^n_s  - \int_0^t  \phi(\mathbb{E}^3 \sigma ^{n_0}(s,\widetilde X_s,  \widetilde \xi_s)) d\widetilde W_s\right\| > c/3\right)
 \\
& \displaystyle +\mathbb P\left(\left\|\int_0^t \phi(\mathbb{E}^3 \sigma^{n_0}(s,\widetilde X_s,  \widetilde \xi_s)) d\widetilde W_s  - \int_0^t  \phi(\mathbb{E}^3 \sigma ^{}(s,\widetilde X_s,  \widetilde \xi_s)) d\widetilde W_s\right\| > c/3\right)
 \\
& \displaystyle =: J^1 + J^2 + J^3.
\end{align*}
Using the assumptions, we apply It\^o-Skorokhod  inequality with any $\delta>0$ and get that
\begin{align*}
& \displaystyle J^1 = \mathbb P\left(\left\|\int_0^t (\phi(\mathbb{E}^3 \sigma^n(s,\widetilde X^n_s,  \widetilde \xi^n_s)) - \phi(\mathbb{E}^3 \sigma ^{n_0}(s,\widetilde X^n_s,  \widetilde \xi^n_s))) d\widetilde W^n_s\right\| > c/3\right)
 \\
& \displaystyle \le \mathbb P\left(\int_0^t \left\|\phi(\mathbb{E}^3 \sigma^n(s,\widetilde X^n_s,  \widetilde \xi^n_s)) - \phi(\mathbb{E}^3 \sigma ^{n_0}(s,\widetilde X^n_s,  \widetilde \xi^n_s))\right\|^2 ds > \delta\right)
 \\
& \displaystyle + \frac9{c^2} \mathbb E \left(\delta \wedge \int_0^{T} \left\|\phi(\mathbb{E}^3 \sigma^n(s,\widetilde X^n_s,  \widetilde \xi^n_s)) - \phi(\mathbb{E}^3 \sigma ^{n_0}(s,\widetilde X^n_s,  \widetilde \xi^n_s))\right\|^2 ds\right).
\end{align*}
Here the second term is small if we choose $\delta>0$ small. Let us consider the first term given $\delta>0$. We have that

\begin{align*}
\mathbb P\left(\int_0^t  \left\|\phi(\mathbb{E}^3 \sigma^n(s,\widetilde X^n_s,  \widetilde \xi^n_s)) - \phi(\mathbb{E}^3 \sigma ^{n_0}(s,\widetilde X^n_s,  \widetilde \xi^n_s))\right\|^2 ds > \delta\right)
 \\
\le \mathbb P\left(C \int_0^t
\left\|\mathbb{E}^3 \sigma^n(s,\widetilde X^n_s,  \widetilde \xi^n_s) - \mathbb{E}^3 \sigma ^{n_0}(s,\widetilde X^n_s,  \widetilde \xi^n_s)\right\|^2 ds > \delta\right)
\end{align*}
By  BCM  inequality,

\begin{align*}
\mathbb P\left(C \int_0^t
\left\|\mathbb{E}^3 \sigma^n(s,\widetilde X^n_s,  \widetilde \xi^n_s) - \mathbb{E}^3 \sigma ^{n_0}(s,\widetilde X^n_s,  \widetilde \xi^n_s)\right\|^2 ds > \delta\right)
 \\
\le (\delta/C)^{-1}
 \mathbb E\int_0^t
\left\|\mathbb{E}^3 (\sigma^n(s,\widetilde X^n_s,  \widetilde \xi^n_s) - \sigma ^{n_0}(s,\widetilde X^n_s,  \widetilde \xi^n_s)\right\|^2 ds
 \\
\le (\delta/C)^{-1}
 \mathbb E\int_0^t
\mathbb{E}^3 \left\|\sigma^n(s,\widetilde X^n_s,  \widetilde \xi^n_s) - \sigma ^{n_0}(s,\widetilde X^n_s,  \widetilde \xi^n_s)\right\|^2 ds
 \\
= (\delta/C)^{-1}
 \mathbb E \mathbb{E}^3\int_0^t
 \left\|\sigma^n(s,\widetilde X^n_s,  \widetilde \xi^n_s) - \sigma ^{n_0}(s,\widetilde X^n_s,  \widetilde \xi^n_s)\right\|^2 ds.
\end{align*}
Convergence of the latter term to zero as $n,n_0\to\infty$ follows from the same considerations as for the drift in the previous step of the proof for the analogous term $I^1$ via Krylov's bound. So, we have
\begin{equation*}
 0\le \;\lim_{n,n_0\to\infty} J^1 \le
\lim_{n,n_0\to\infty} \mathbb E \mathbb{E}^3\int_0^t
 \left\|\sigma^n(s,\widetilde X^n_s,  \widetilde \xi^n_s) - \sigma ^{n_0}(s,\widetilde X^n_s,  \widetilde \xi^n_s)\right\|^2 ds = 0.
\end{equation*}

The term $J^2$ converges to zero by Skorokhod's Lemma \ref{app1}:
\begin{align*}
\int_0^t \phi(\mathbb{E}^3 \sigma^{n_0}(s,\widetilde X^n_s,  \widetilde \xi^n_s)) d\widetilde W^n_s  \stackrel{\mathbb P}{\to} \int_0^t  \phi(\mathbb{E}^3 \sigma ^{n_0}(s,\widetilde X_s,  \widetilde \xi_s)) d\widetilde W_s, \quad n\to\infty,
\end{align*}
with
$
f^n_s:= \phi(\mathbb{E}^3 \sigma^{n_0}(s,\widetilde  X^n_s, \widetilde  \xi^n_s)), \quad f^0_s:= \phi(\mathbb{E}^3 \sigma^{n_0}(s,\widetilde  X_s, \widetilde  \xi_s))
$
in this lemma.

Consider $J^3$:
\begin{align*}
J^3 = \mathbb P\left(\left\|\int_0^t (\phi(\mathbb{E}^3 \sigma^{n_0}(s,\widetilde X_s,  \widetilde \xi_s)) - \phi(\mathbb{E}^3 \sigma ^{}(s,\widetilde X_s,  \widetilde \xi_s))) d\widetilde W_s\right\| > c/3\right)
 \\
\le \mathbb P\left(\int_0^t \left\|\phi(\mathbb{E}^3 \sigma^{n_0}(s,\widetilde X_s,  \widetilde \xi_s)) - \phi(\mathbb{E}^3 \sigma ^{}(s,\widetilde X_s,  \widetilde \xi_s))\right\|^2 ds > \delta\right)
 \\
+ \frac9{c^2} \mathbb E \left(\delta \wedge \int_0^{T} \left\|\phi(\mathbb{E}^3 \sigma^n(s,\widetilde X_s,  \widetilde \xi_s)) - \phi(\mathbb{E}^3 \sigma ^{n_0}(s,\widetilde X_s,  \widetilde \xi_s))\right\|^2 ds\right).
\end{align*}
Similarly to $J^1$, the second  term in the last sum  is small for small $\delta$. For the first one we have, similarly to $J^1$,
\begin{align*}
\mathbb P\left(\int_0^t \left\|\phi(\mathbb{E}^3 \sigma^{n_0}(s,\widetilde X_s,  \widetilde \xi_s)) - \phi(\mathbb{E}^3 \sigma ^{}(s,\widetilde X_s,  \widetilde \xi_s))\right\|^2 ds > \delta\right)
 \\
\le \mathbb P\left(C \int_0^t
\left\|\mathbb{E}^3 \sigma^{n_0}(s,\widetilde X_s,  \widetilde \xi_s) - \mathbb{E}^3 \sigma ^{}(s,\widetilde X_s,  \widetilde \xi_s)\right\|^2 ds > \delta\right)
 \\
\le (\delta/C)^{-1}
\mathbb E \mathbb{E}^3\int_0^t
 \left\|\sigma^{n_0}(s,\widetilde X_s,  \widetilde \xi_s) - \sigma ^{}(s,\widetilde X_s,  \widetilde \xi_s)\right\|^2 ds.
\end{align*}
Convergence of this term to zero as $n_0\to\infty$ follows from  Lemma   \ref{lekrybd}, similarly to the analogous convergence of $I^3$ in the previous step. So,
$$
\lim_{n_0\to\infty}J^3 = 0.
$$
This finishes the proof of the desired bound (\ref{sest0}). Thus, a weak solution of the equation (\ref{e1})--(\ref{e200}) exists in the case of $d_1=d$ and under the  assumption (\ref{si1}) instead of (\ref{si}). For bounded coefficients and under (\ref{phiA0})--(\ref{psiA0}) the Proposition \ref{pro22} is proved.

\medskip

\noindent
{\bf 6}.
Now consider the general case of unbounded coefficients satisfying the linear growth condition in $x$ along with (\ref{fifi})--(\ref{psipsi}).
Let
\begin{align*}
& \displaystyle \sigma^n(t,x,y) = (\sigma(t,x,y) \mathbf{1}(|x|\le n)
+ \sigma(0,0,0) \mathbf{1}(|x| > n))*\phi_n(x),
\\
& \displaystyle b^n(t,x,y) = b(t,x,y) \mathbf{1}(|x|\le n).
\end{align*}
Note that the function $\sigma^n$ is Borel measurable,  bounded, uniformly non-degenerate, and smooth (at least, Lipschitz) in $x$. Denote by $X^n$ a  solution of the equation
$$
X^n_t = X_0 + \int_0^t \overline{B}^n[s,X^n_s,\mu^n_s]\,ds + \int_0^t  \overline{\Sigma}^n[s,X^n_s,\mu^n_s]\,dW^n_s,
$$
or
$$
X^n_t = X_0 + \int_0^t \psi(\mathbb E^3 b^n[s,X^n_s,\xi^n_s])\,ds + \int_0^t \phi(\mathbb E^3 \sigma^n[s,X^n_s,\xi^n_s])\,dW^n_s,
$$
with
\begin{align*}
  \overline{\Sigma}^n[s,x,\mu] = \phi\left(\int \sigma^n(s,x,y)\,\mu(dy)\right)
= \phi(\Sigma^n[s,x,\mu]),
 \\
\overline{ B}^n[s,x,\mu] = \psi\left(\int b^n(s,x,y)\,\mu(dy)\right)
= \phi(B^n[s,x,\mu]),
\end{align*}
and where $(\xi^n)$ are independent of $(X^n,W^n)$ processes with the same distributions as $X^n$ on some independent probability space. (Note that we use the same notations as in the step 1; however, the approximations $b^n$ and $\sigma^n$ here are different; so we define them as the new ones.) This weak (in fact, strong) solution exists for each $n$   due to the Proposition \ref{pro22} for bounded coefficients. We will use again Skorokhod's technique; note that we could not apply it in one step because in our method of justifying the previous steps the boundedness of both coefficients is essential.

A priori moment inequalities of Lemma \ref{lebds} hold true uniformly with respect to $n$. Hence, by Skorokhod's Lemma \ref{app2}, choose a subsequence of equivalently distributed triples $(\widetilde X^{n_k}, \widetilde \xi^{n_k}, \widetilde W^{n_k})$ converging in probability for any $s$ to some limiting triple $(\widetilde X^{}, \widetilde \xi^{}, \widetilde W^{})$. Here $\widetilde W^{}$ is a Wiener process of dimension $d_1$ (cf. with \cite[Chapter 2]{Kry}, where, however, dimensions are equal, but it does not affect the conclusion). Convergence in the
equation
$$
\widetilde X^{n_k}_t = \widetilde X_0 + \int_0^t \psi(\mathbb E^3 b^{n_k}[s,\widetilde X^{n_k}_s,\widetilde \xi^{n_k}_s])\,ds + \int_0^t \phi(\mathbb E^3 \sigma^{n_k}[s,\widetilde X^{n_k}_s,\widetilde \xi^{n_k}_s])\,d\widetilde W^{n_k}_s
$$
towards the limiting equation
$$
\widetilde X^{}_t = \widetilde x_0 + \int_0^t \psi(\mathbb E^3 B[s,\widetilde X^{}_s,\widetilde \xi^{}_s])\,ds + \int_0^t \phi(\mathbb E^3 \sigma^{}[s,\widetilde X^{}_s,\widetilde \xi^{}_s])\,d\widetilde W_s
$$
follows from the same calculus as in the proof of the previous steps of Proposition \ref{pro22} with the only difference that now $\sigma$ may also be unbounded; however, this does require some additional care, in particular,  because we want to use Krylov's bounds stated for unbounded coefficients (and do not forget about $\phi$ and $\psi$). Yet,  in fact, we can use Krylov's bound  from Lemma \ref{lekrybd}.

First of all, recall that a priori bounds
 \eqref{kol1} and \eqref{kol2} hold true with constants  not depending on $n$.
Now, by Skorokhod's Lemma   \ref{app2},  on some probability space we have some equivalent processes $(\widetilde X^{n'}_t,\widetilde \xi^{n'}_t, \widetilde W^{n'}_t)$ and a limiting triple $(\widetilde X^{}_t,\widetilde \xi^{}_t, \widetilde W^{}_t)$ such that for any $t$,
\begin{equation*}
(\widetilde X^{n'}_t,\widetilde \xi^{n'}_t, \widetilde W^{n'}_t) \stackrel{\mathbb P}{\to} (\widetilde X^{}_t,\widetilde  \xi^{}_t, \widetilde W^{}_t).
\end{equation*}
By virtue of   a priori bounds for \(\widetilde W^n\), the process \(\widetilde W \) is continuous and it is a Wiener process.
Also, the limits are adapted to the corresponding filtration $\widetilde {\cal F}_t:= \bigvee_{n}{\cal F}^{(n)}_t$, and \(\widetilde W\) is  a Wiener process with respect to this filtration. In particular,  related  Lebesgue and stochastic integrals are all well defined.
Moreover, by virtue of the uniform estimates (\ref{kol2}), the limit \((\widetilde X^{}_t,\widetilde  \xi^{}_t)\) may be also regarded as  continuous due to Kolmogorov's continuity theorem  because   a priori bounds  (\ref{kol1}) -- (\ref{kol2}) remain valid for the limiting processes \(\widetilde X, \widetilde \xi\). In particular, it is useful to note for the sequel that
\begin{eqnarray*}
\sup_{0\le t\le T} \mathbb{E} |\widetilde X_t|^2 \le C_T (1+\mathbb E |x_0|^2)\; \text{and}\; \mathbb{E} \sup_{0\le t\le T} |\widetilde X_t|^2 \le C_T (1+\mathbb E |x_0|^2).
\end{eqnarray*}

\medskip

\noindent
{\bf 7}.
\noindent
Let us now show that
\begin{equation}\label{limNb}
\int_0^t \psi(\mathbb{E}^3 b^{n'}(s,\widetilde X^{n'}_s, \widetilde \xi^{n'}_s))ds  \stackrel{\mathbb P}{\to} \int_0^t \psi( \mathbb{E}^3 b(s,\widetilde X_s, \widetilde \xi_s))ds,
\end{equation}
and
\begin{equation}\label{limNs}
\int_0^t \phi(\mathbb{E}^3 \sigma^{n'}(s,\widetilde X^{n'}_s,  \widetilde \xi^{n'}_s)) d\widetilde W^{n'}_s  \stackrel{\mathbb P}{\to} \int_0^t  \phi(\mathbb{E}^3 \sigma(s,\widetilde X_s, \widetilde  \xi_s)) d\widetilde W_s, \quad n'\to\infty.
\end{equation}

\medskip

Denote $\widetilde{R}=R-1.$ Given any $\epsilon >0$, and simplifying notation by renaming $n'\rightarrow n$, for any $t\le T$ by virtue of the BCM  inequality we conclude that for any $\epsilon>0$ there exists $  R$ such that
\[
\mathbb E \mathbf{1}(\gamma_{n,  \widetilde R} \wedge \gamma_{  \widetilde R}  \le  T) < \epsilon.
\]
Further at one place we will need more precise estimates:
\begin{equation*}\label{precise1}
\mathbb P(\gamma^X_{n,\widetilde R}\le T) \le
   (\widetilde R)^{-2} {\mathbb E \sup_{0\le t\le T} |\widetilde X^n_t|^2}
\le (\widetilde R)^{-2}{C(1+\mathbb E |x_0|^2)},
\end{equation*}
and
\begin{equation*}\label{precise2}
\mathbb P(\gamma^X_{\widetilde R}\le T) \le
(\widetilde R)^{-2}{\mathbb E \sup_{0\le t\le T} |\widetilde X_t|^2}
\le (\widetilde R)^{-2}{C(1+\mathbb E |x_0|^2)},
\end{equation*}
by virtue of the BCM  inequality. Hence,
\begin{align*}
\mathbb P(\gamma^X_{n,\widetilde R}\wedge \gamma^X_{\widetilde R}\le T) \le
2C(\widetilde R)^{-2}{(1+\mathbb E |x_0|^2)}.
\end{align*}
Now, we estimate, with $\widetilde c = c/C$, $m=m_2$,
\begin{align*}
& \displaystyle \mathbb P\left(\left|\int_0^t \psi(\mathbb{E}^3 b^{n}(s,\widetilde X^{n}_s, \widetilde \xi^{n}_s))\,ds -  \int_0^t \psi(\mathbb{E}^3 b(s,\widetilde X_s, \widetilde \xi_s))\,ds\right|>c\right)
\\
& \displaystyle  \le
\mathbb P\left(C \int_0^t (1+(|\widetilde X^{n}_s|\vee |\widetilde X^{}_s|)^m)\left|\mathbb{E}^3 \left(b^{n}(s,\widetilde X^{n}_s, \widetilde \xi^{n}_s) -  b(s,\widetilde X_s, \widetilde \xi_s)\right)\right|
\,ds > c\right)
\\
& \displaystyle  =
\mathbb P\left(\int_0^t (1+(|\widetilde X^{n}_s|\vee |\widetilde X^{}_s|)^m)\mathbb{E}^3\!\left|\left(\mathbf{1}(\gamma_{n,\widetilde R} \wedge \gamma_{\widetilde R}   \le  T)\!+\!\mathbf{1}(\gamma_{n,\widetilde R} \wedge \gamma_{\widetilde R}   >  T) \right)
 \right. \right.\\ \left.\left.\right.\right.\\
& \displaystyle \hspace{4cm}
\left. \left.\times
\left(b^{n}(s,\widetilde X^{n}_s, \widetilde \xi^{n}_s) -  b(s,\widetilde X_s, \widetilde \xi_s)\right)\right|ds>\widetilde c\right)
\\
& \displaystyle
\le \mathbb P\left(\int_0^t(1+(|\widetilde X^{n}_s|\vee |\widetilde X^{}_s|)^m) \mathbb{E}^3\!\mathbf{1}(\gamma_{n,\widetilde R} \wedge \gamma_{\widetilde R}   \le  T)\left|\left(b^{n}(s,\widetilde X^{n}_s, \widetilde \xi^{n}_s) -  b(s,\widetilde X_s, \widetilde \xi_s)\right)\right|ds>\widetilde c/2\right)
\\
& \displaystyle
+\mathbb P\left(\int_0^t(1+(|\widetilde X^{n}_s|\vee |\widetilde X^{}_s|)^m)\mathbb{E}^3\!\mathbf{1}(\gamma_{n,\widetilde R} \wedge \gamma_{\widetilde R}   >  T) \left|\left(b^{n}(s,\widetilde X^{n}_s, \widetilde \xi^{n}_s) -  b(s,\widetilde X_s, \widetilde \xi_s)\right)\right|ds>\widetilde c/2\right)
\\
& \displaystyle
\le \mathbb P\left(\mathbf{1}(\gamma^X_{n,\widetilde R} \wedge \gamma^X_{\widetilde R} \le  T) \int_0^t (1+(|\widetilde X^{n}_s|\vee |\widetilde X^{}_s|)^m)\mathbb{E}^3\!
\left|\left(b^{n}(s,\widetilde X^{n}_s, \widetilde \xi^{n}_s) -  b(s,\widetilde X_s, \widetilde \xi_s)\right)\right|ds>\widetilde c/2\right)
\\
& \displaystyle + \mathbb P\left(
\int_0^t (1+(|\widetilde X^{n}_s|\vee |\widetilde X^{}_s|)^m)\mathbb{E}^3\!\mathbf{1}(\gamma^\xi_{n,\widetilde R} \wedge \gamma^\xi_{\widetilde R}   \le  T)
\left|\left(b^{n}(s,\widetilde X^{n}_s, \widetilde \xi^{n}_s) -  b(s,\widetilde X_s, \widetilde \xi_s)\right)\right|ds>\widetilde c/2\right)
\\
& \displaystyle
+\mathbb P\left(\mathbf{1}(\gamma^X_{n,\widetilde R} \wedge \gamma^X_{\widetilde R}>T)\int_0^t(1+(|\widetilde X^{n}_s|\vee |\widetilde X^{}_s|)^m)\mathbb{E}^3\!\mathbf{1}(\gamma^\xi_{n,\widetilde R} \wedge \gamma^\xi_{\widetilde R}   >  T)
 \right. \\ \left.\right.\\
& \displaystyle \hspace{4cm}
\left. \times \left|\left(b^{n}(s,\widetilde X^{n}_s, \widetilde \xi^{n}_s) -  b(s,\widetilde X_s, \widetilde \xi_s)\right)\right|ds>\widetilde c/2\right)
 \\
& \displaystyle \le \mathbb P(\gamma^X_{n,\widetilde R} \wedge \gamma^X_{\widetilde R}   \le  T)
\\
& \displaystyle + \mathbb P\left(
\int_0^t (1+(|\widetilde X^{n}_s|\vee |\widetilde X^{}_s|)^m)\mathbb{E}^3\!\mathbf{1}(\gamma^\xi_{n,\widetilde R} \wedge \gamma^\xi_{\widetilde R}   \le  T)
\left|\left(b^{n}(s,\widetilde X^{n}_s, \widetilde \xi^{n}_s) -  b(s,\widetilde X_s, \widetilde \xi_s)\right)\right|ds>\widetilde c/2\right)
\\
 & \displaystyle
+ \mathbb P \left(\gamma^X_{n,\widetilde R} \wedge \gamma^X_{\widetilde R}>T;  \int_0^t (1+(|\widetilde X^{n}_s|\vee |\widetilde X^{}_s|)^m)\mathbb{E}^3 \mathbf{1}(\gamma^\xi_{n,\widetilde R} \wedge \gamma^\xi_{\widetilde R}  >  T)
 \right. \\ \left.\right.\\
& \displaystyle \hspace{4cm}
\left. \times
 \left|\left(b^{n}(s,\widetilde X^{n}_s, \widetilde \xi^{n}_s)  - b(s,\widetilde X_s,\widetilde  \xi_s)\right)\right|ds > \frac{\widetilde c}{2}\right) =: L^1 + L^2 + L^3.
\end{align*}
Here the first term $L^1$ does not exceed $\epsilon$ if $R$ is large enough, uniformly with respect to $n$.

Consider the term $L^2$. We estimate it using the linear growth in $x$ assumption:
\begin{eqnarray*}
& \displaystyle \mathbb P\left(
\int_0^t (1+(|\widetilde X^{n}_s|\vee |\widetilde X^{}_s|)^m)\mathbb{E}^3\!\mathbf{1}(\gamma^\xi_{n,\widetilde R} \wedge \gamma^\xi_{\widetilde R}   \le  T)
\left|\left(b^{n}(s,\widetilde X^{n}_s, \widetilde \xi^{n}_s) -  b(s,\widetilde X_s, \widetilde \xi_s)\right)\right|ds>\widetilde c/2\right)
\\
& \displaystyle
\le  \mathbb P\left(\mathbb{E}^3\!\mathbf{1}(\gamma^\xi_{n,\widetilde R} \wedge \gamma^\xi_{\widetilde R}   \le  T)
\int_0^t (1+(|\widetilde X^{n}_s|\vee |\widetilde X^{}_s|)^{m+1})
ds>\widetilde c/2\right).
\end{eqnarray*}
Since the processes $\widetilde X^{n}_s$ together with their limit $\widetilde X^{}_s$ are bounded in probability uniformly in $n$, and because
$$
\mathbb{E}^3\!\mathbf{1}(\gamma^\xi_{n,\widetilde R} \wedge \gamma^\xi_{\widetilde R}   \le  T) \to 0, \quad R\to\infty,
$$
we conclude that
$$
\lim_{R \to \infty}\sup_n L^2 = 0.
$$

Further, consider $L^3$. Let us fix some  $n_0$ and let $n>n_0$. Replacing $\widetilde c/2$ by $c$ for simplicity, we have for any $t\le T$,
\begin{eqnarray*}\label{3terms}
& \displaystyle \mathbb P \left(\gamma^X_{n,\widetilde R} \wedge \gamma^X_{\widetilde R}  >  T; \, \int_0^t (1+(|\widetilde X^{n}_s|\vee |\widetilde X^{}_s|)^m) \mathbb{E}^3 \mathbf{1}(\gamma^\xi_{n,\widetilde R} \wedge \gamma^\xi_{\widetilde R}  >  T)
 \right. \\ \left.\right.\\
& \displaystyle \hspace{4cm}
\left. \times \left|\left(b^{n}(s,\widetilde X^{n}_s, \widetilde \xi^{n}_s)  - b(s,\widetilde X_s,\widetilde  \xi_s)\right)\right|ds > c\right)
 \nonumber \\\nonumber \\\nonumber
& \displaystyle \le \mathbb P \left(\gamma^X_{n,\widetilde R} \wedge \gamma^X_{\widetilde R}  >  T;  \int_0^t (1+(|\widetilde X^{n}_s|\vee |\widetilde X^{}_s|)^m)\mathbb{E}^3 \mathbf{1}(\gamma^\xi_{n,\widetilde R} \wedge \gamma^\xi_{\widetilde R }  >  T)
 \right. \\ \left.\right.\\
& \displaystyle \hspace{4cm}
\left. \times \left|\left(b^{n}(s,\widetilde X^{n}_s, \widetilde \xi^{n}_s)  - b^{n_0}(s,\widetilde X^n_s,\widetilde  \xi^n_s)\right)\right|ds > \frac{c}3\right)
 \nonumber \\\nonumber \\\nonumber
& \displaystyle + \mathbb P \left(\gamma^X_{n,\widetilde R} \wedge \gamma^X_{\widetilde R}  >  T;  \int_0^t (1+(|\widetilde X^{n}_s|\vee |\widetilde X^{}_s|)^m) \mathbb{E}^3 \mathbf{1}(\gamma^\xi_{n,\widetilde R} \wedge \gamma^\xi_{\widetilde R}  >  T)
 \right. \\ \left.\right.\\
& \displaystyle \hspace{4cm}
\left. \times
 \left|\left(b^{n_0}(s,\widetilde X^{n}_s, \widetilde \xi^{n}_s)  - b^{n_0}(s,\widetilde X_s,\widetilde  \xi_s)\right)\right|ds > \frac{c}3\right)
 \nonumber \\\nonumber \\\nonumber
& \displaystyle + \mathbb P \left(\gamma^X_{n,\widetilde R} \wedge \gamma^X_{\widetilde R}  >  T;  \int_0^t (1+(|\widetilde X^{n}_s|\vee |\widetilde X^{}_s|)^m)\mathbb{E}^3 \mathbf{1}(\gamma^\xi_{n,\widetilde R} \wedge \gamma^\xi_{\widetilde R}  >  T)
 \right. \\ \left.\right.\\
& \displaystyle \hspace{4cm}
\left. \times  \left|\left(b^{n_0}(s,\widetilde X^{}_s, \widetilde \xi^{}_s)  - b(s,\widetilde X_s,\widetilde  \xi_s)\right)\right|ds > \frac{c}3\right)
  \\ \nonumber\\
& \displaystyle =: M^1 + M^2 + M^3.
\end{eqnarray*}
Denote as earlier
\[
g^{n,n_0}(s,x,\xi):=b^{n}(s,x,\xi)- b^{n_0}(s,x,\xi), \quad
g^{n_0}(s,x,\xi):=b^{n_0}(s,x,\xi)- b^{}(s,x,\xi).
\]
Let $\alpha\le 1/(2m+2)$. Note that
$$
\mathbb E
(1+(\sup_s|\widetilde X^{n}_s|\vee \sup_s|\widetilde X^{}_s|)^{2m\alpha})
\le C < \infty,
$$
and that this constant $C$ is uniform in $n$ and does not depend on $R$.
The first summand $M^1$  may be estimated by the BCM inequality  as
\begin{eqnarray}\label{i1}
M^1 \le C\, \mathbb E \left( \mathbf{1}(\gamma^X_{n,\widetilde R} > T, \gamma^X_{\widetilde R}> T)\int_0^T(1+(|\widetilde X^{n}_s|\vee |\widetilde X^{}_s|)^m) \right.
  \nonumber \\  \left. \right. \nonumber \\ \left.
\times \mathbb{E}^3 \mathbf{1}(\gamma^\xi_{n,\widetilde R} \wedge \gamma^\xi_{\widetilde R}  >  T)|b^{n}(s,\widetilde X^{n}_s,\widetilde  \xi^{n}_s)- b^{n_0}(s,\widetilde X^n_s, \widetilde \xi^n_s)|\,ds\right)^\alpha
  \nonumber \\ \nonumber\\ \nonumber
\le C\, \mathbb E \mathbf{1}(\gamma^X_{n,\widetilde R} > T, \gamma^X_{\widetilde R}> T)
(1+(\sup_s|\widetilde X^{n}_s|\vee \sup_s|\widetilde X^{}_s|)^{m\alpha})
 \nonumber \\  \nonumber \\  \nonumber
\times  \left( \int_0^T \mathbb{E}^3 \mathbf{1}(\gamma^\xi_{n,\widetilde R} \wedge \gamma^\xi_{\widetilde R}  >  T)|b^{n}(s,\widetilde X^{n}_s,\widetilde  \xi^{n}_s)- b^{n_0}(s,\widetilde X^n_s, \widetilde \xi^n_s)|\,ds\right)^\alpha
  \nonumber \\ \nonumber\\ \nonumber
\stackrel{\mbox{CBS}}{\le} C\, \left(\mathbb E \mathbf{1}(\gamma^X_{n,\widetilde R} > T, \gamma^X_{\widetilde R}> T)
(1+(\sup_s|\widetilde X^{n}_s|\vee \sup_s|\widetilde X^{}_s|)^{2m\alpha})\right)^{1/2}
 \nonumber \\  \nonumber \\  \nonumber
\times  \left(\mathbb E \left( \int_0^T \mathbb{E}^3 \mathbf{1}(\gamma^\xi_{n,\widetilde R} \wedge \gamma^\xi_{\widetilde R}  >  T)|b^{n}(s,\widetilde X^{n}_s,\widetilde  \xi^{n}_s)- b^{n_0}(s,\widetilde X^n_s, \widetilde \xi^n_s)|\,ds\right)^{2\alpha}\right)^{1/2}
 \nonumber \\ \nonumber\\ \nonumber
\stackrel{\mbox{H\"older}}{\le} C\, \left( \int_0^T \mathbb{E} \mathbf{1}(\gamma^\xi_{n,\widetilde R} \wedge \gamma^\xi_{\widetilde R}  >  T)|b^{n}(s,\widetilde X^{n}_s,\widetilde  \xi^{n}_s)- b^{n_0}(s,\widetilde X^n_s, \widetilde \xi^n_s)|\,ds\right)^{\alpha}
  \nonumber \\ 
  \\ \nonumber
=  C\, \left( \int_0^T \mathbb{E} \mathbf{1}(\gamma^\xi_{n,\widetilde R} \wedge \gamma^\xi_{\widetilde R}  >  T)|g^{n,n_0}(s,\widetilde X^{n}_{s\wedge \gamma_{n,\widetilde R}},\widetilde  \xi^{n}_{s\wedge \gamma_{n,\widetilde R}})|\,ds\right)^{\alpha}.
\end{eqnarray}
Here the couple $(\widetilde X^{n}_{s\wedge \gamma_{n,\widetilde R}},\widetilde  \xi^{n}_{s\wedge \gamma_{n,\widetilde R}})$ is a stopped diffusion with coefficients bounded by norm in state variables $(x,\xi)$ by the value $C \widetilde R$ uniformly with respect to $n$, and with the diffusion coefficient  uniformly non-degenerate.  

Denote by $\widehat B^{n,\widetilde R} [s,x,\mu]$ and $\widehat \Sigma^{n,\widetilde R} [s,x,\mu]$  bounded vector and matrix functions in $x$ respectively, with $\widehat \Sigma^{n,\widetilde R} [s,x,\mu]$ uniformly non-degenerate and smooth (e.g., \(C^1\)), such that
\[
\widehat B^{n,\widetilde R} [s,x,\mu] = B^n [s,x,\mu], \;\; \widehat \Sigma^{n,\widetilde R}[s,x,\mu] = \Sigma^n[s,x,\mu], \quad |x|\le \widetilde R.
\]
Let
\((\widehat X^{n}_{s}) = (\widehat X^{n,\widetilde R}_{s})\)
be a  strong  solution of the It\^o equation,
\begin{equation*}
d\widehat X^n_t = \psi(\widehat B^{n,\widetilde R}[t,\widehat X^n_t, \mu^n_t])\,dt
+ \phi(\widehat \Sigma^{n,\widetilde R} [t,\widehat X^n_t, \mu^n_t])\,d\widetilde W^{n}_t, \quad \widehat X^n_0=x_0,
\end{equation*}
where $\mu^n_t$ is still the marginal distribution of \(X^n_t\) and \(\widetilde X^n_t\). Let also $\widehat \xi^n$ be an equivalent independent copy of the process \(\widehat X^n_t\) satisfying the equation
\begin{equation*}
d\widehat \xi^n_t = \psi(\widehat B^{n,\widetilde R}[t,\widehat \xi^n_t, \mu_t])dt + \phi(\widehat \Sigma^{n,\widetilde R} [t,\widehat\xi^n_t, \mu_t]d\widetilde W^{',n}_t), \;\; t\ge 0, \quad
{\cal L}(\xi^n_0)={\cal L}(x_0).
\end{equation*}
We may assume that the Wiener processes $\widetilde W^{',\;n}$ here are the same as in the equation (\ref{exi0}) for $\widetilde \xi^n$. (Emphasize that solutions $\widehat\xi^n$ are strong ones; this is why we have mollified $\sigma$ in the variable $x$).
Then it follows that on $[0,t\wedge \gamma_{n,\widetilde R}]$ the processes $\widetilde X^n$ and $\widehat X^n$ coincide, see
\cite[Theorem 6.2.1(v)]{Kry-ln}, as well as $\widetilde \xi^n$ coincide with $\widehat \xi^n$.
 We highlight that the same stopping times $\gamma_{n,\widetilde R}$ can be used. Then the bound  for $I^1$ in the second line of (\ref{i1}) may be rewritten as
\begin{eqnarray}\label{i11}
& \displaystyle M^1 \le C\, \left(\mathbb E  \mathbf{1}(\gamma_{n,\widetilde R} > T, \gamma_{\widetilde R}>  T)\int_0^T  |g^{n,n_0}(s,\widehat  X^{n}_{s\wedge \gamma_{n,\widetilde R}},\widehat  \xi^{n}_{s\wedge \gamma_{n,\widetilde R}})| \,ds\right)^\alpha
 \nonumber\\ \nonumber
& \displaystyle \le  C\, \left(\mathbb E \int_0^T  |g^{n,n_0}(s,\widehat  X^{n}_{s\wedge \gamma_{n,\widetilde R}},\widehat  \xi^{n}_{s\wedge \gamma_{n,\widetilde R}})| \,ds\right)^\alpha.
\end{eqnarray}
The values of the function $g^{n,n_0}(s,x,\xi)$ outside the set $\{(x,\xi):\, (|x|\vee |\xi|)\le \widetilde R\}$
are not relevant to the evaluation of the expression in the second line of (\ref{i11}). So, without loss of generality we may assume for our goal that $g^{n,n_0}(s,x,\xi)$ vanishes outside of this ball. Then, by  Krylov's estimate from Lemma \ref{lekrybd}  we obtain with some constant $N$ depending on the dimension $d$ and on $R$ through the ellipticity constants of the diffusion matrix and the sup-norm of the drift with some $N_R$,
\begin{eqnarray*}
& \displaystyle M^1 \le C\, \left(\mathbb E  \int_0^T  |g^{n,n_0}(s,\widehat  X^{n}_{s\wedge \gamma_{n,\widetilde R}},\widehat  \xi^{n}_{s\wedge \gamma_{n,\widetilde R}})| \,ds\right)^\alpha
\\
& \displaystyle =  N_R\, \left(\int_0^T \int_{|x|\le \widetilde R}\int_{|\xi|\le \widetilde R} |b^{n}(s,x,\xi) - b^{n_0}(s,x,\xi)|^{2d+1}\,dx d\xi ds\right)^{\frac{\alpha}{2d+1}}
\\
& \displaystyle \le
N_R\, \left(\int_0^T \int_{|x|\le \widetilde R}\int_{|\xi|\le \widetilde R} |b^{n}(s,x,\xi) - b^{}(s,x,\xi)|^{2d+1}\,dx d\xi ds\right)^{\frac{\alpha}{2d+1}}
 \\
& \displaystyle + N_R\, \left(\int_0^T \int_{|x|\le \widetilde R}\int_{|\xi|\le \widetilde R} |b^{n_0}(s,x,\xi) - b^{}(s,x,\xi)|^{2d+1}\,dx d\xi ds\right)^{\frac{\alpha}{2d+1}}
\to 0, \quad n, n_0 \to \infty,
\end{eqnarray*}
for each $R$, by virtue of the properties of mollified functions. Hence, for any $R$,
$$
\lim_{n,n_0\to \infty}M^1 = 0.
$$

\medskip

Let
 $0\le w(x, \xi) \le 1$ be any continuous function which equals $1$ for every $|x|\vee |\xi|\le \widetilde R$ and zero for every $|x|\vee |\xi|>R$, and $g^{n_0}_R(s,x,\xi) := g^{n_0}(s,x,\xi) w(x,\xi)$.

Further, the second term $M^2$  for any \(0\le t\le T\) admits the following upper bound:
\begin{eqnarray*}
&\displaystyle \!\!\! M^2 \!\le\!
\mathbb P \!\left(\gamma^X_{n,\widetilde R} \!\wedge \!\gamma^X_{\widetilde R}  >  T; \!\int_0^t (1+(|\widehat X^{n}_{s\wedge \gamma_{n,\widetilde R}}|\vee |\widehat X^{}_{s\wedge \gamma_{n,\widetilde R}}|)^m) \right.
 \\\left.\right.\\
&\left.
\times
 \mathbb{E}^3 \mathbf{1}(\gamma^\xi_{n,\widetilde R} \! \wedge\! \gamma^\xi_{\widetilde R}\! >\! T)\!\left|\left(b^{n_0}(s,\widetilde X^{n}_s, \widetilde \xi^{n}_s)ds \! -\! b^{n_0}(s,\widetilde X_s,\widetilde  \xi_s)\right)\right|\!ds \!>\! \frac{c}3\right)
\\
&\displaystyle \le \frac3{c}\,(1+\widetilde R^m) \mathbb E \mathbb{E}^3 \mathbf{1}(\gamma_{n,\widetilde R}> T, \gamma_{\widetilde R} > T)\int_0^T  |b^{n_0}(s,\widetilde X^{n}_s, \widetilde \xi^{n}_s)- b^{n_0}(s,\widetilde X_s, \widetilde \xi_s)|\,ds
 \\
&\displaystyle \le \frac3{c}\,(1+\widetilde R^m) \mathbb E \mathbb{E}^3 \mathbf{1}(\gamma_{n,\widetilde R}> T, \gamma_{\widetilde R} > T)
 \\
&\displaystyle \times
\int_0^T  |b^{n_0}(s,\widetilde X^{n}_s, \widetilde \xi^{n}_s) w(\widetilde X^{n}_s, \widetilde \xi^{n}_s) - b^{n_0}(s,\widetilde X_s, \widetilde \xi_s) w(\widetilde X^{}_s, \widetilde \xi^{}_s)|\,ds.
\end{eqnarray*}
Let us show that the latter expression tends to zero as $n\to \infty$.

The random variable
\(\displaystyle \int_0^T  |b^{n_0}(s,\widetilde X^{n}_s, \widetilde \xi^{n}_s) w(\widetilde X^{n}_s, \widetilde \xi^{n}_s) - b^{n_0}(s,\widetilde X_s, \widetilde \xi_s) w(\widetilde X^{}_s, \widetilde \xi^{}_s)|\,ds\) is bounded uniformly in $n$.

So, it suffices to show that
\begin{equation}\label{why}
\int_0^T    |b^{n_0}(s,\widetilde X^{n}_s, \widetilde \xi^{n}_s) w(\widetilde X^{n}_s, \widetilde \xi^{n}_s) - b^{n_0}(s,\widetilde X_s, \widetilde \xi_s) w(\widetilde X^{}_s, \widetilde \xi^{}_s)|\,ds \to 0, \quad n\to \infty,
\end{equation}
in probability. Indeed, (\ref{why}) follows from the standard trick in Krylov's estimates: firstly we approximate again the bounded function $b^{n_0}(s,x,\xi)w(x,\xi)$
by smooth (continuous suffices) in $(x,\xi)$ function $\bar b^{n_0}(s,x,\xi)$ also with a compact support, so that their difference in the $L_{2d+1}$ norm is small. Then
$$
\int_0^T    |\bar b^{n_0}(s,\widetilde X^{n}_s, \widetilde \xi^{n}_s)  - \bar b^{n_0}(s,\widetilde X_s, \widetilde \xi_s) |\,ds \to 0, \quad n\to \infty,
$$
due to the continuity of $\bar b^{n_0}$.
Also, the difference
$$
\mathbb E\int_0^T    |b^{n_0}(s,\widetilde X^{n}_s, \widetilde \xi^{n}_s) w(\widetilde X^{n}_s, \widetilde \xi^{n}_s) - \bar b^{n_0}(s,\widetilde X^{n}_s, \widetilde \xi^{n}_s)  |\,ds
$$
is small due to  Krylov's estimate (Lemma \ref{lekrybd}); it does not exceed some constant $N_R$ multiplied by $\|b^{n_0}w - \bar b^{n_0}\|_{[0,T] \times B_R\times B_R}$.

For any Borel function $g_R$ Krylov's bound (Lemma \ref{lekrybd}) implies the inequality
\begin{align*}
\sup_n \mathbb E
\int_0^T  |g_R(s,\widetilde X^{n}_s, \widetilde  \xi^{n}_s)|\,ds
\le N_R \|g_R\|_{L_{2d+1}([0,T]\times B_R \times B_R)},
\end{align*}
where, as before,  $B_R=\{z\in \mathbb{R}^d: |z|\le R\}$.
Assume that $g_R$ is continuous in $(x,y)$ and vanishing outside $B_{R}\times B_{R}$.
Then in the limit we obtain  by Fatou lemma  that
\begin{equation}\label{kbunbd}
\mathbb E
\int_0^T  |g_R^{}(s,\widetilde X^{}_s, \widetilde  \xi^{}_s)|\,ds \le N_R \|g_R\|_{L_{2d+1}([0,T]\times B_R \times B_R)}.
\end{equation}
 
This inequality for unbounded coefficients was already justified in the comments after (\ref{kb}).


Therefore,
$$
\lim_{n,n_0\to\infty}M^2 = 0.
$$

\medskip

Let us consider the third term $M^3$. We have that
\begin{equation*}
\begin{split}
M^3 &= \mathbb P \!\left(\!\gamma^X_{n,\widetilde R} \!\wedge \!\gamma^X_{\widetilde R}\!  >\!  T;  \!\int_0^T \! (1+(|\widetilde X^{n}_s|\vee |\widetilde X^{}_s|)^m) \mathbb{E}^3\! \mathbf{1}(\gamma^\xi_{n,\widetilde R}\! \wedge \!\gamma^\xi_{\widetilde R} \! > \! T)\!
 \right.
\\
&\times \left.\left|\left(b^{n_0}(s,\widetilde X^{}_s, \widetilde \xi^{}_s) \! - \! b^{}(s,\widetilde X_s,\widetilde  \xi_s)\right)\right|\!ds\! > \!
\frac{c}3\right)
\\
&\le C(1+\widetilde R^m)\,\mathbb E
\mathbf{1}(\gamma_{\widetilde R}\wedge \gamma_{n,\widetilde R} > T)
\int_0^T  |g^{n_0}(s,\widetilde X^{}_s, \widetilde  \xi^{}_s)|\,ds
\\
&\le C(1+\widetilde R^m)\,\mathbb E
\int_0^T  |g_R^{n_0}(s,\widetilde X^{}_s, \widetilde  \xi^{}_s)|\,ds.
\end{split}
\end{equation*}
Now, by the same reasons as for the term $M^2$, we obtain
\begin{equation*}
\begin{split}
M^3&\le C(1+\widetilde R^m)\,
\mathbb E
\mathbf{1}(\gamma_{\widetilde R}\wedge \gamma_{n,\widetilde R} > T)
\int_0^T  |g^{n_0}(s,\widetilde X^{}_s, \widetilde  \xi^{}_s)|\,ds
\\
&\le C(1+\widetilde R^m)\,
\mathbb E
\int_0^T  |g^{n_0}_R(s,\widetilde X^{}_s, \widetilde  \xi^{}_s)|\,ds
\\
&\le C(1+\widetilde R^m)\,N_R \|g^{n_0}_R\|_{L_{2d+1}([0,T]\times B_R \times B_R)}
\end{split}
\end{equation*}
 due to (\ref{kbunbd}). Hence,
$$
\lim_{n_0\to\infty} M^3 = 0.
$$
The convergence  (\ref{limNb}) is, thus, proved.

\medskip

\noindent
{\bf 8}. Now let us consider convergence of stochastic integrals in (\ref{limNs}). Our goal is an estimate the next probability similarly to   the drift and Lebesgue integrals above:
\begin{equation}\label{sest}
\mathbb P\left(\left|\int_0^t \phi(\mathbb{E}^3 \sigma^n(s,\widetilde X^n_s,  \widetilde \xi^n_s)) d\widetilde W^n_s  - \int_0^t  \phi(\mathbb{E}^3 \sigma(s,\widetilde X_s,  \widetilde \xi_s)) d\widetilde W_s\right| > c\right)
 < C \epsilon,
\end{equation}
for any $c, \epsilon>0$ and $n$ large enough. As we know from the step 5 in the proof of Proposition \ref{pro22}, the task is similar to the convergence of Lebesgue integrals, and  we mainly show how to tackle the additional obstacle due to different Wiener processes \( \widetilde W \) and \(\widetilde W^n\) in the stochastic integrals. A tool for this  is Skorokhod's Lemma  \ref{app1}, as before. However, it is not applicable directly because our processes may be unbounded, so we should overcome this with the help of a truncation  which reduces the problem to bounded processes.

By virtue of \cite[Theorem 6.2.1(v)]{Kry-ln} and similarly to the calculus for Lebesgue integrals in the previous steps, yet using second moments instead of the first ones by the evident reason we estimate,

\begin{align*}
\mathbb P\left(\left|\int_0^t \phi(\mathbb{E}^3 \sigma^n(s,\widetilde X^n_s,  \widetilde \xi^n_s)) d\widetilde W^n_s  - \int_0^t  \phi(\mathbb{E}^3 \sigma(s,\widetilde X_s,  \widetilde \xi_s)) d\widetilde W_s\right| > c\right)
\\
= \mathbb P\left(\gamma^X_{\widetilde R}\wedge \gamma^X_{n,\widetilde R} \le T; \left|\int_0^t \phi(\mathbb{E}^3 \sigma^n(s,\widetilde X^n_s,  \widetilde \xi^n_s)) d\widetilde W^n_s  - \int_0^t  \phi(\mathbb{E}^3 \sigma(s,\widetilde X_s,  \widetilde \xi_s)) d\widetilde W_s\right| > c\right)
\\
+\mathbb P\left(\gamma^X_{\widetilde R}\wedge \gamma^X_{n,\widetilde R} > T; \left|\int_0^t \phi(\mathbb{E}^3 \sigma^n(s,\widetilde X^n_s,  \widetilde \xi^n_s)) d\widetilde W^n_s  - \int_0^t  \phi(\mathbb{E}^3 \sigma(s,\widetilde X_s,  \widetilde \xi_s)) d\widetilde W_s\right| > c\right)
\\
\le \mathbb P(\gamma^X_{\widetilde R}\wedge \gamma^X_{n,\widetilde R} \le T)
\\
+\mathbb P\left(\gamma^X_{\widetilde R}\wedge \gamma^X_{n,\widetilde R} > T; \left|\int_0^t \phi(\mathbb{E}^3 \sigma^n(s,\widehat X^n_s,  \widehat \xi^n_s)) d\widetilde W^n_s  - \int_0^t  \phi(\mathbb{E}^3 \sigma(s,\widetilde X_s, \widetilde \xi_s)) d\widetilde W_s\right| > c\right).
\end{align*}
The first term here $\mathbb P(\gamma^X_{\widetilde R}\wedge \gamma^X_{n,\widetilde R} \le  T)$ is small if $R$ is large enough, as we have seen in the earlier steps. It remains to consider the second term. We will evaluate it as follows:

\begin{align*}
& \displaystyle \mathbb P(\gamma^X_{\widetilde R}\wedge \gamma^X_{n,\widetilde R} > T; \left|\int_0^t \phi(\mathbb{E}^3 \sigma^n(s,\widehat X^n_s,  \widehat \xi^n_s)) d\widetilde W^n_s  - \int_0^t  \phi(\mathbb{E}^3 \sigma(s,\widetilde X_s,  \widetilde \xi_s)) d\widetilde W_s\right| > c)
\\
& \displaystyle \le \mathbb P\left(\gamma^X_{\widetilde R}\wedge \gamma^X_{n,\widetilde R} > T; \left|\int_0^T \phi(\mathbb{E}^3 \sigma^n(s,\widehat X^n_s,  \widehat \xi^n_s)) d\widetilde W^n_s  - \int_0^t  \phi(\mathbb{E}^3 \sigma^{n_0}(s,\widehat X^n_s,  \widehat \xi^n_s)) d\widetilde W^n_s\right| > c/3\right)
\\
& \displaystyle + \mathbb P\left(\gamma^X_{\widetilde R}\wedge \gamma^X_{n,\widetilde R} > T; \left|\int_0^T \phi(\mathbb{E}^3 \sigma^{n_0}(s,\widetilde X^n_s,  \widetilde \xi^n_s)) d\widetilde W^n_s  - \int_0^t  \phi(\mathbb{E}^3 \sigma^{n_0} (s,\widetilde X_s,  \widetilde \xi_s)) d\widetilde W_s\right| > c/3\right)
\\
& \displaystyle + \mathbb P\left(\gamma^X_{\widetilde R}\wedge \gamma^X_{n,\widetilde R} > T; \left|\int_0^T \phi(\mathbb{E}^3 \sigma^{n_0}(s,\widetilde X_s,  \widetilde \xi_s)) d\widetilde W_s  - \int_0^t  \phi(\mathbb{E}^3 \sigma^{} (s,\widetilde X_s,  \widetilde \xi_s)) d\widetilde W_s\right| > c/3\right)
\\
& \displaystyle =: J^1 + J^2 + J^3.
\end{align*}
Now for all three terms the evaluation is similar to that for the drift term, except for one difference. We start with $J^1$:
\begin{align*}
J^1
=\mathbb P\left(\gamma^X_{\widetilde R}\wedge \gamma^X_{n,\widetilde R} > T; \left|\int_0^T (\phi(\mathbb{E}^3 \sigma^n(s,\widehat X^n_s,  \widetilde \xi^n_s)) d\widetilde W^n_s  - \phi(\mathbb{E}^3 \sigma^{n_0}(s,\widehat X^n_s,  \widetilde \xi^n_s))) d\widetilde W^n_s\right| > c/3\right)
\\
\le \mathbb P\left(\left|\int_0^{T\wedge \gamma^X_{\widetilde R}\wedge \gamma^X_{n,\widetilde R}} (\phi(\mathbb{E}^3 \sigma^n(s,\widehat X^n_s,  \widetilde \xi^n_s))  - \phi(\mathbb{E}^3 \sigma^{n_0}(s,\widehat X^n_s,  \widetilde \xi^n_s))) d\widetilde W^n_s\right| > c/3\right)
\\
=\mathbb P\left(\left|\int_0^{T} \mathbf{1}(s <\gamma^X_{\widetilde R}\wedge \gamma^X_{n,\widetilde R}) (\phi(\mathbb{E}^3 \sigma^n(s,\widehat X^n_s,  \widetilde \xi^n_s))   - \phi(\mathbb{E}^3 \sigma^{n_0}(s,\widehat X^n_s,  \widetilde \xi^n_s))) d\widetilde W^n_s\right| > c/3\right).
\end{align*}
By virtue of the  simplified version of the It\^o-Skorokhod's inequality (see \cite[Theorem 6.3.5]{Kry-ln}, we get that  for any $\delta>0$
\begin{align*}
& \displaystyle \mathbb P\left(\left|\int_0^{T} \mathbf{1}(s <\gamma^X_{\widetilde R}\wedge \gamma^X_{n,\widetilde R}) (\phi(\mathbb{E}^3 \sigma^n(s,\widehat X^n_s,  \widetilde \xi^n_s))   - \phi(\mathbb{E}^3 \sigma^{n_0}(s,\widehat X^n_s,  \widetilde \xi^n_s))) d\widetilde W^n_s\right| > c/3\right)
\\
& \displaystyle \le \mathbb P\left(\int_0^{T} \mathbf{1}(s <\gamma^X_{\widetilde R}\wedge \gamma^X_{n,\widetilde R}) \left\|\phi(\mathbb{E}^3 \sigma^n(s,\widehat X^n_s,  \widetilde \xi^n_s))  - \phi(\mathbb{E}^3 \sigma^{n_0}(s,\widehat X^n_s,  \widetilde \xi^n_s))\right\|^2 ds > \delta\right)
\\
& \displaystyle + \frac9{c^2} \mathbb E \left(\delta \wedge \int_0^{T} \mathbf{1}(s <\gamma^X_{\widetilde R}\wedge \gamma^X_{n,\widetilde R})\left\|\phi(\mathbb{E}^3 \sigma^n(s,\widehat X^n_s,  \widetilde \xi^n_s))  - \phi(\mathbb{E}^3 \sigma^{n_0}(s,\widehat X^n_s,  \widetilde \xi^n_s))\right\|^2 ds\right)
\\
& \displaystyle =: F^1 + F^2.
\end{align*}
The term $F^2$ is small if we choose $\delta$ small enough. Let us consider $F^1$. It can be bounded as follows:
\begin{gather*}
\mathbb P\left(\int_0^{T} \mathbf{1}(s <\gamma^X_{\widetilde R}\wedge \gamma^X_{n,\widetilde R}) \left\|\phi(\mathbb{E}^3 \sigma^n(s,\widehat X^n_s,  \widetilde \xi^n_s))  - \phi(\mathbb{E}^3 \sigma^{n_0}(s,\widehat X^n_s,  \widetilde \xi^n_s))\right\|^2 ds > \delta\right)\\ \le
\mathbb P\left(\int_0^{T} \mathbf{1}(s <\gamma^X_{\widetilde R}\wedge \gamma^X_{n,\widetilde R}) (1+|\widehat X^n_s|^m)^2 \left\|\mathbb{E}^3 (\sigma^n(s,\widehat X^n_s,  \widetilde \xi^n_s)  -  \sigma^{n_0}(s,\widehat X^n_s,  \widetilde \xi^n_s))\right\|^2 ds > \delta\right)
 \\
=
\mathbb P\bigg(\int_0^{T} \mathbf{1}(s <\gamma^X_{\widetilde R}\wedge \gamma^X_{n,\widetilde R}) (1+|\widehat X^n_s|^m)^2 \left\|\mathbb{E}^3
(\mathbf{1}(T \ge \gamma^\xi_{\widetilde R}\wedge \gamma^\xi_{n,\widetilde R})
 \right. \\ \left.
+ \mathbf{1}(T < \gamma^\xi_{\widetilde R}\wedge \gamma^\xi_{n,\widetilde R}))
\times (\sigma^n(s,\widehat X^n_s,  \widetilde \xi^n_s)  -
\sigma^{n_0}(s,\widehat X^n_s,  \widetilde \xi^n_s))\right\|^2 ds > \delta\bigg)\\
\le
\mathbb P\bigg(\int_0^{T} \mathbf{1}(s <\gamma^X_{\widetilde R}\wedge \gamma^X_{n,\widetilde R})\mathbb{E}^3  (1+|\widehat X^n_s|^m)^2
\mathbf{1}(T \ge \gamma^\xi_{\widetilde R}\wedge \gamma^\xi_{n,\widetilde R})\\
\times \|\sigma^n(s,\widehat X^n_s,  \widetilde \xi^n_s)  -
\sigma^{n_0}(s,\widehat X^n_s,  \widetilde \xi^n_s)||^2 ds > \delta/2\bigg)
 \\
+ \mathbb P\bigg(\int_0^{T} \mathbf{1}(s <\gamma^X_{\widetilde R}\wedge \gamma^X_{n,\widetilde R}) (1+|\widehat X^n_s|^m)^2 \mathbb{E}^3
 \mathbf{1}(T < \gamma^\xi_{\widetilde R}\wedge \gamma^\xi_{n,\widetilde R})
 \\
 \times \|\sigma^n(s,\widehat X^n_s,  \widetilde \xi^n_s)  -
\sigma^{n_0}(s,\widehat X^n_s,  \widetilde \xi^n_s)\||^2 ds > \delta/2\bigg)
 =: S^1 + S^2.
\end{gather*}
For $S^2$ we can replace $\widetilde\xi$ by $\widehat \xi$ and then use BCM inequality as in the drift evaluation, to obtain an $L_{4d+1}$-bound (not $L_{2d+1}$ because of $\|\sigma^n - \sigma^{n_0}\|^2$) by virtue of Krylov's bound (see Lemma \ref{lekrybd}):
\begin{gather*}
S^2 = \mathbb P\bigg(\int_0^{T} \mathbf{1}(s <\gamma^X_{\widetilde R}\wedge \gamma^X_{n,\widetilde R}) (1+|\widehat X^n_s|^m)^2 \mathbb{E}^3
 \mathbf{1}(T < \gamma^\xi_{\widetilde R}\wedge \gamma^\xi_{n,\widetilde R})\\
 \times \|\sigma^n(s,\widehat X^n_s,  \widehat \xi^n_s)  -
\sigma^{n_0}(s,\widehat X^n_s,  \widehat \xi^n_s)\|^2 ds > \delta/2\bigg)\\
  \le C_R \|g^{n,n_0}\|^2_{L_{4d+1}([0,T]\times B_R\times B_R)} \to 0, \quad n,n_0 \to\infty.
\end{gather*}
For $S^1$ let $\beta = 1/(4m+4)$ and also use the BCM inequality:
\begin{gather*}
\mathbb P\bigg(\int_0^{T} \mathbf{1}(s <\gamma^X_{\widetilde R}\wedge \gamma^X_{n,\widetilde R})\mathbb{E}^3  (1+|\widehat X^n_s|^m)^2
\mathbf{1}(T \ge \gamma^\xi_{\widetilde R}\wedge \gamma^\xi_{n,\widetilde R})
 \times \|\sigma^n(s,\widehat X^n_s,  \widetilde \xi^n_s)\\  -
\sigma^{n_0}(s,\widehat X^n_s,  \widetilde \xi^n_s)||^2 ds > \delta/2\bigg)
 \le C \mathbb E \bigg(\int_0^{T} \mathbf{1}(s <\gamma^X_{\widetilde R}\wedge \gamma^X_{n,\widetilde R})\mathbb{E}^3  (1\\+|\widehat X^n_s|^m)^2
\mathbf{1}(T \ge \gamma^\xi_{\widetilde R}\wedge \gamma^\xi_{n,\widetilde R})
  \times \|\sigma^n(s,\widehat X^n_s,  \widetilde \xi^n_s)  -
\sigma^{n_0}(s,\widehat X^n_s,  \widetilde \xi^n_s)||^2 ds\bigg)^\beta
 \\
\stackrel{\mbox{CBS}}{\le}
C (T \mathbb E \sup_s (1+|\widehat X^n_s|)^{})^{1/2} (\mathbb P (T \ge \gamma^\xi_{\widetilde R}\wedge \gamma^\xi_{n,\widetilde R}))^{\beta/2}.
\end{gather*}
Due to   a priori bounds from Lemma \ref{lebds}, the first multiplier here is bounded, while the second is small if $R$ is large enough, all uniformly in $n$. Hence
 $
\lim_{R\to\infty} S^1 = 0.
$
Therefore,
$$
\lim_{n, n_0\to \infty}J^1 = 0.
$$
 Both other two terms $J^2$ and $J^3$ are tackled similarly to the drift terms in the previous step with the use of (\ref{kbunbd}) and It\^o's isometry.

This proves the desired bound (\ref{sest}).
So, (\ref{limNb}) and (\ref{limNs}) hold true, and thus, the proposition is proved.

\subsection{Proof of Theorem \ref{thm1}}

We will use
a hint from \cite[section 4]{K-V} in order to reduce the statement to the case considered earlier in the Proposition \ref{pro22}. However, for the reader's convenience we repeat the details.

Denote $\widetilde \Sigma[t,x,\mu] := \sqrt{A[t,x,\mu]}$, where
$A[t,x,\mu]:= \Sigma[t,x,\mu]\Sigma^*[t,x,\mu]$,
and   suppose  that there exists a  weak  solution $\widetilde X$ of the equation
\begin{equation*}
\widetilde X_t = x_0 + \int_0^t B[s,\widetilde X_s, \mu_s]ds +
\int_0^t  \widetilde \Sigma[s,\widetilde X_s,  \mu_s] d\widetilde W_s,
\end{equation*}
with some $d$-dimensional Wiener process $(\widetilde W_t, \, t\ge 0)$ on some
probability space and where $\mu_s$ stands for the distribution of $\widetilde
X_s$.

Existence of this weak solution is already justified in the Proposition \ref{pro22}.
 Here the positive-definite matrix-function $\sqrt{A[t,x,\mu]}$ is defined via the Cauchy-Riesz-Dunford  formula (\ref{contour}), or, equivalently, (\ref{sqrta2}).

\medskip

Further, without loss of generality we may and will assume that on the same
probability space there exists another \emph{independent} $d_1$-dimensional
Wiener process $(\overline{W}_t, \, t\ge 0)$. Let $I$ denote a $d_1\times
d_1$-dimensional unit matrix and let
\begin{equation*}
p[s,x,\mu] := \widetilde \Sigma[s,x,\mu]^{-1}\, \Sigma[s,x,\mu].
\end{equation*}
Note that the matrix $\widetilde \Sigma[s,x,\mu]$ is symmetric and that
\begin{align*}
p^*p[s,x,\mu] = \Sigma^*[s,x,\mu]  (\widetilde \Sigma^*[s,x,\mu])^{-1}\widetilde \Sigma[s,x,\mu]^{-1}\, \Sigma[s,x,\mu]
\\
= \Sigma^*[s,x,\mu]  (A)^{-1}[s,x,\mu]  \Sigma[s,x,\mu] ,
\\
p^*[s,x,\mu] p[s,x,\mu] p^*[s,x,\mu] p[s,x,\mu]
\\
= \Sigma^*[s,x,\mu]  (A)^{-1}[s,x,\mu]  \Sigma[s,x,\mu]  \Sigma^*[s,x,\mu]  (A)^{-1}[s,x,\mu]  \Sigma[s,x,\mu]
\\
= \Sigma^* (A)^{-1} (A)(A)^{-1} \Sigma[s,x,\mu] = \Sigma^* (A)^{-1} \Sigma[s,x,\mu],
\end{align*}
and let
$$
W^0_t:= \int_0^t p^*[s,\widetilde X_s, \mu_s] \, d\widetilde W_s
+ \int_0^t (I - p^*[s,\widetilde X_s, \mu_s] p[s,\widetilde X_s, \mu_s])
\, d\overline{W}_s.
$$
Notice that
\begin{eqnarray*}
\Sigma[s,x,\mu] p^*[s,x,\mu] = A[s,x,\mu] (A[s,x,\mu])^{-1/2} =(A[s,x,\mu])^{1/2},
\\
\Sigma[s,x,\mu] p^*[s,x,\mu]p[s,x,\mu] = (A[s,x,\mu])^{1/2}p[s,x,\mu]
\\
= (A[s,x,\mu])^{1/2}(a[s,x,\mu])^{-1/2} \Sigma[s,x,\mu] = \Sigma[s,x,\mu].
\end{eqnarray*}
Due to the
multivariate L\'evy characterization theorem
this implies that $W^0$ is a $d_1$-dimensional Wiener
process, since its matrix angle characteristic (also known as a matrix angle bracket) equals
\begin{eqnarray*}
\langle W^0, W^0\rangle_t = \int_0^t p^*p[s,\widetilde X_s,\mu_s]\,ds + \int_0^t  (I-p^* p[s,\widetilde X_s,\mu_s])^*(I-p^* p[s,\widetilde X_s,\mu_s]) \,ds
\\
= \int_0^t  (p^* p[s,\widetilde X_s,\mu_s] + I - 2 p^* p[s,\widetilde X_s,\mu_s] + p^* pp^* p[s,\widetilde X_s,\mu_s])\,ds
\\
= \int_0^t  (I - p^* p[s,\widetilde X_s,\mu_s] + p^* pp^* p[s,\widetilde X_s,\mu_s])\,ds
\\
= \int_0^t  (I - \Sigma^* (A)^{-1} \Sigma[s,\widetilde X_s,\mu_s] + \Sigma^* (A)^{-1}(A)(A)^{-1} \Sigma[s,\widetilde X_s,\mu_s])\,ds= \int_0^t I \,ds = t\,I.
\end{eqnarray*}
Next, due to the stochastic integration rules (see \cite{IkedaWatanabe}),
\begin{eqnarray}\label{last1}
 \int_0^t \Sigma[s,\widetilde X_s, \mu_s]\,d W^0_s = \int_0^t  \Sigma p^*[s,\widetilde X_s, \mu_s] \, d\widetilde W_s + \int_0^t  \Sigma (I-p^* p)[s,\widetilde X_s, \mu_s]\, d\overline{W}_s
 \nonumber\\ \nonumber
= \int_0^t  (A)^{1/2}[s,\widetilde X_s, \mu_s]\, d\widetilde W_s
= \int_0^t  \widetilde \Sigma[s,\widetilde X_s, \mu_s] \, d\widetilde W_s  = \widetilde X_t - x_0 - \int_0^t B[s,\widetilde X_s, \mu_s]\,ds.
\end{eqnarray}
In other words, \((\widetilde X, W^0)\) is a  weak  solution of the equation (\ref{e1}).
It remains to notice that since we did not change measures, $\mu_s$ is still the distribution of $\widetilde X_s$
by the assumption.
The proof of the Theorem \ref{thm1} is thus completed.

\begin{remark}
There is a non-rigorous view that  for SDE solutions everything related to weak solutions and weak uniqueness depends only on the matrix $\sigma^*\sigma$ and not on $\sigma$ itself. This is not precise. Firstly, for strong solutions this is not true because regularity such as Lipschitz condition or even a simple continuity may fail for a badly chosen square root, let us forget about non-Borel square roots. Secondly, even for weak solutions in the absence of non-degeneracy and if the square root is not continuous, there is no guarantee that weak solution exists for   any  square root. Also, existing results about weak solutions and weak uniqueness, see \cite{Chiang, Funaki}, impose conditions on \(\sigma\) and not on \(\sigma\sigma^*\).  Hence, a vague  ``common knowledge'' is not sufficient and had to show the calculus.
\end{remark}

\section{Strong solutions; strong and weak uniqueness}\label{sec:se}
\subsection{On strong existence}
In this section it is shown that strong solution of the equation  (\ref{e1})--(\ref{e200}) exists
under appropriate conditions. Emphasize that we do
not claim strong   uniqueness  in this section, but only strong existence in the sense of the Definition {\ref{Def1}}. We also notice for
interested readers that in \cite{Ver82} the assumption of continuity in time was dropped in comparison to
\cite{Ver80}; so, just a certain  local Lipschitz condition suffices for our aim.

\begin{Proposition}\label{thm4}
Let \(\mathbb E |x_0|^4 < \infty\).
Also, let the coefficients $b$ and
$\sigma$ satisfy all conditions of the Theorem \ref{thm1} and the non-degeneracy assumption (\ref{si1}),  and let just $\sigma$ be Lipschitz in $x$ uniformly with
respect to $s$ and locally with respect to $y$,
\begin{equation}\label{loclips}
\|\sigma(t,x,y) - \sigma(t,x',y)\| \le C(1+|y|^2)|x-x'|.
\end{equation}

Then
the equation
(\ref{e1})--(\ref{e200}) has a strong solution and, moreover, every
solution is strong and, in particular, solution may be constructed
on any probability space equipped with a $d_1$-dimensional
Wiener process.
\end{Proposition}

This result is likely to be a common knowledge. A brief sketch of the proof is presented below for completeness and because the authors were unable to find an exact reference.

\medskip

\noindent
{\bf 1.}
First of all, note that   weak solutions exist  and    a priori bounds (\ref{kol1})--(\ref{kol2}) are valid.
Considerations are based on the results from \cite{Ver80} and \cite{Ver82} about strong solutions for SDEs
for a Borel measurable drift which is assumed bounded or with a linear growth in both papers. Since weak
solution does exist, whatever is its distribution $\mu$, the process $X$ may be considered
as an ordinary SDE with coefficients depending on time,
$$
\widetilde b(t,x) = B[t,x,\mu_t], \qquad \widetilde \sigma(t,x) =
\Sigma[t,x,\mu_t],
$$
and, hence,
\begin{equation}\label{e22}
dX_t = \widetilde b(t,X_t) dt + \widetilde \sigma(t,X_t)
dW_t, \qquad
X_0=x.
\end{equation}
Recall that according to   Lemma \ref{cor1}, the new coefficients $\widetilde b(t,x)$ and $\widetilde \sigma(t,x)$ are Borel measurable.

\medskip

\noindent
{\bf 2.}
Now in order to establish strong existence it suffices to verify that the new coefficients $\widetilde b$ and
$\widetilde \sigma$ both satisfy linear growth in $x$ condition uniform in time, and $\widetilde \sigma$ is Lipschitz
continuous in $x$ and remains uniformly non-degenerate.

\medskip

\noindent
We have, for any $T>0$ and $0\le t \le T$,
\begin{eqnarray*}
|\widetilde b(t,x)|
=  \left|\int b(t,x,y)\,
\mu_t(dy)\right|
 \le C  \int (1 + |x|) \, \mu_t(dy)   =  C\,(1 +
|x|).
\end{eqnarray*}
Similarly, it also follows that
\begin{eqnarray*}
\|\widetilde \sigma(t,x)\|
\le C \int (1 + |x|) \, \mu_t(dy)   = C\,(1 + |x|).
\end{eqnarray*}
Further,   by virtue of \eqref{loclips},
\begin{eqnarray*}
\|\widetilde \sigma(t,x) - \widetilde \sigma(t,x')\| =  \|\Sigma[t,x,
\mu_t] - \Sigma[t,x',\mu_t]\|
 \\
= \left\|\int \sigma(t,x,y)\, \mu_t(dy)  - \int \sigma(t,x',y)\,
\mu_t(dy) \right\|
 \\
\le C \,|x-x'|\, \int (1 + |y|^2) \, \mu_t(dy)  \le C_{T}\,
\,|x-x'|.
\end{eqnarray*}
The uniform non-degeneracy of $\sigma$, hence, also of $\sigma\sigma^*$, follows from the inequality (\ref{si1}) by integration with respect to $\mu_t$.

These properties suffice for the local  strong uniqueness of solution of the equation (\ref{e1})--(\ref{e200}) by  virtue of the results in \cite{Ver80}. However, because weak solution is well-defined for all values of time, strong uniqueness is global.
According to the Yamada-Watanabe principle (\cite{YamadaWatanabe}), any solution of the equation (\ref{e200}) is strong. So, any solution of the original equation (\ref{e1})  is also strong.
This completes the proof of the Proposition~\ref{thm4}.

\begin{remark}\label{rem2}
Notice that as a solution of the ``linearized''
equation (\ref{e22}), $X$ is pathwise unique, but
so far it is not known if this implies the same property
for $X$ as a solution of (\ref{e1}), unless {\em weak uniqueness} for the equation
(\ref{e1}) has been established. In a restricted framework this will be done in Theorem  \ref{thm5} below.
\end{remark}

\begin{remark}
In the case of $d=1$, Lipschitz condition may
be relaxed to H\"older of order $1/2$ and, actually, a little bit further by using techniques from \cite{YamadaWatanabe} and \cite{Ver-1d}.
\end{remark}

\subsection{Weak uniqueness}
In this section weak uniqueness will be shown for the equation (\ref{e1}) -- (\ref{e200}) under appropriate
conditions. This result  requires only a Borel  measurability
of the drift with respect to the state variable $x$, although, it
assumes that diffusion $\sigma$ does not depend on $y$
along with some additional continuity condition in $x$ and  non-degeneracy. The drift may be unbounded in the state variable $x$.

\begin{theorem}\label{thm5a}
	Let
	\(\mathbb E \exp(r |x_0|^2)<\infty\) for some \(r>0\), and let the functions
	$b$ and $\sigma$ be Borel
	measurable,  and
	$$
	\sigma(s,x,y) \equiv \sigma (s,x),
	$$
	that is, \(\sigma\) does not depend on the variable \(y\);
	let $\sigma$ satisfy the non-degeneracy assumption (\ref{si1});
	let $d_1 = d$, the matrix $\sigma$ be bounded, symmetric and invertible, and let there exist $C>0$ such that the function
	$$
	\widetilde
	B[s,x,\mu] := \sigma^{-1}(s,x)\, B[s,x,\mu]
	$$
	satisfies the linear growth condition: there is $C>0$ such that for all $x\in \mathbb{R}^d$,
	\begin{equation}\label{bgrow}
	\sup_{s,\mu}|\widetilde B[s,x,\mu]| \le C(1+|x|).
	\end{equation}
	Also assume that the matrix-function $\sigma(t,x)$ satisfies the uniform continuity condition in $x$  which guarantees that the equation
	\begin{equation}\label{itosde}
	dX^0_t = \sigma(t,X^0_t)\,dW_t, \quad X^0_0 = \xi,
	\end{equation}
	with an ${\cal F}_0$-measurable initial data $\xi$ possessing a given distribution $\mu_0 = {\cal L}(x_0)$,
	has a unique weak solution  (see  \cite{Ver80, Ver82}).
	 Then under the assumptions of Theorem \ref{thm1}  solution of the equation (\ref{e1})--(\ref{e200}) is weakly unique.
\end{theorem}

\begin{remark}\label{rem3}
In case of $d=1$ continuity of $\sigma$ in $x$ is not needed.
	Under the additional assumption of boundedness of \(\widetilde b\) or \(\widetilde B\) an exponential moment of the initial value \(x_0\) is not necessary and can be replaced by the fourth moment,
or even weaker.
\end{remark}

\subsection{Proof of Theorem \ref{thm5a}}
Denote by $X^0_t$ any (weakly unique) weak  solution of the It\^o equation~(\ref{itosde}).
Note that
$$
dW^{}_t = \sigma^{-1}(t,X_t^0)dX^0_t.
$$

\noindent
{\bf 1.}
Recall that under the
assumptions of the theorem, any solution of the equation (\ref{e1})--(\ref{e200}) is
strong by virtue of the Proposition \ref{thm4}.
Hence, it suffices to
show weak uniqueness, after which strong uniqueness for this equation will follow from
strong uniqueness for the ``linearized'' equation (\ref{e22}). We will show this weak uniqueness  by contradiction. Suppose there
are two solutions $(X^1,W^1)$ and $(X^2,W^2)$ of the equation (\ref{e1}) with distributions $\mu^1$
and $\mu^2$ respectively in the space of trajectories
$C[0,\infty; \mathbb R^d]$:
\begin{equation}\label{new1}
dX^1_t = \sigma(t, X^1_t)\,d W^1_t
+B[t,X^1_t, \mu^1_t]\,dt, \quad X^1_0=\xi^1,
\end{equation}
and
\begin{equation}\label{new2}
dX^2_t = \sigma(t, X^2_t)\,dW^2_t
+B[t,X^2_t, \mu^2_t]\,dt, \quad X^2_0=\xi^2,
\end{equation}
respectively, with ${\cal L}(\xi^1) = {\cal L}(\xi^2)= {\cal L}(x_0)$.
Yet, under the present setting it will be shown that firstly $\mu^1 = \mu^2$ and secondly $X^1 = X^2$ a.s.
Note that both $X^1$ and $X^2$ are Markov processes (\cite{Kry_selection}).

\medskip

Both solutions $(X^i,\mu^i)$   in the weak sense  may be obtained from the Wiener process $W$ and solution $X^0$ of the equation without the drift (\ref{itosde}) via
Girsanov's transformations using the following  stochastic exponents:
$$
\gamma^i_T = \exp\left(\int_0^T \widetilde B[s,X^0_s,\mu^i_s]\,dW_s -
\frac12
\int_0^T |\widetilde B[s,X^0_s,\mu^i_s]|^2 \,ds\right), \quad i=1,2,
$$
where $\widetilde
b(t,x,y) := \sigma^{-1}(t,x)\, b(t,x,y)$, $\widetilde
B[t,x,\mu] := \sigma^{-1}(t,x)\, B[t,x,\mu]$, $|\widetilde B|$ stands for the modulus of the vector $\widetilde B$,
 and $\widetilde B[s,X^0_s,\mu^i_s]\,dW_s$ is understood as a scalar product, \(\sum_{j=1}^{d}\widetilde B^j[s,X^0_s,\mu^i_s] d\widetilde W^j_s\).

\medskip

It is well-known that in the case of bounded $\widetilde B$ the random variables $\gamma^i_T$, $i=1,2$, are probability densities due to Girsanov  theorem (see, e.g., \cite[Theorem 6.8.8]{Kry-ln}).
So, till the step 4 we assume \(\widetilde B\) bounded; note that in this case we have,

\begin{equation} \label{liptvr}
|\widetilde B[s,x,\mu] - \widetilde B[s,x,\nu]| \le C \|\mu-\nu\|_{TV}.
\end{equation}
Indeed,
\begin{align*}
|\widetilde B[s,x,\mu] - \widetilde B[s,x,\nu]|
= |\sigma^{-1}(s,x)\int (b(s,x,y)\mu(dy) - b(s,x,y)\nu(dy)])|
 \\\\
= |\sigma^{-1}(s,x)\int (b(s,x,y)(\mu(dy) - \nu(dy))])|
\le \int |\sigma^{-1}(s,x)b(s,x,y)|\,|\mu(dy) - \nu(dy)|
 \\\\
\le  \|\sigma^{-1}(s,x)b(s,x,\cdot)\|_B \|\mu-\nu\|_{TV}.
\end{align*}

The calculus with a bounded  $\widetilde B$  is needed so as to explain the idea which will be further expanded to the case without this restriction. Also this will justify the statement in the Remark \ref{rem3}.

\medskip

Denote
$$
\widetilde W^1_t := W_t - \int_0^t \widetilde B[s, X^0_s, \mu_s^{1}]\,ds, \quad 0\le t\le T.
$$
This is a new Wiener process on $[0,T]$ under the probability measure $P^{\gamma^1}$ defined by its density as $(dP^{\gamma^1}/dP)(\omega) = \gamma^1_T$.
Then, on the same interval $[0,T]$, on the probability space with a Wiener process $(\Omega, {\cal F}, (\widetilde W^1_t, F^{}_t), \mathbb P^{\gamma^1})$, the process $(X^0_t, \, 0\le t\le T)$ satisfies the equation,
\begin{equation}\label{x01}
\begin{split}
dX^0_t &= \sigma(t, X^0_t)d \widetilde W^1_t
+ \sigma(t, X^0_t) \widetilde B[t,X^0_t, \mu^1_t]dt
\\
&= \sigma(t, X^0_t)d \widetilde W^1_t
+B[t,X^0_t, \mu^1_t]\,dt, \hspace{1cm}
\end{split}
\end{equation}
with the initial condition $X^0_0 =x_0$. In other words, the process $X^0$ on $[0,T]$ satisfies the equation (\ref{new1}), just with another Wiener process and under another probability measure. However, given $\mu^1_t, \, 0\le t\le T$, this solution considered as a solution of It\^o's  or ``linearized''  equation is weakly unique.
This is a well-known fact for bounded coefficients due to the results on uniqueness for solutions of parabolic equations, see \cite{Lad}. For unbounded coefficients under the linear growth conditions this follows by truncation and via stopping times in a standard way. Further, this uniqueness for $X^0$ implies weak uniqueness for the pair $(X^0,W)$, see \cite{Bass} et al.

So, the pair $(X^0_t, \widetilde W^1_t, \, 0\le t\le T)$ has the same distribution under the measure $\mathbb P^{\gamma^1}$ as the pair $(X^1_t, W_t, \, 0\le t\le T)$ under the measure $\mathbb P$. Therefore, the marginal distribution of $X^0_t$ under the measure $\mathbb P^{\gamma^1}$ equals $\mu^1_t$, i.e., the couple $(X^0_t, \mu^1_t)$ under $\mathbb P^{\gamma^1}$ solves the McKean-Vlasov equation (\ref{e1}), that is, it is equivalent to the pair $(X^1_t, \mu^1_t, \, 0\le t\le T)$ under the measure $\mathbb P$.

Note for the sequel that $d\widetilde W^1_t$ admits a representation
$$
d\widetilde W^1_t = \sigma^{-1}(t,X_t^0)\,dX^0_t
- \sigma^{-1}(t,X_t^0) B[t,X^0_t, \mu^1_t]\,dt
= \sigma^{-1}(t,X_t^0)dX^0_t - \widetilde B[t,X^0_t, \mu^1_t]\,dt,
$$
or, equivalently,
$$
\sigma^{-1}(t,X_t^0)dX^0_t = d\widetilde W^1_t + \widetilde B[t,X^0_t, \mu^1_t]\,dt.
$$
Similarly, let
$$
\widetilde W^2_t := W_t - \int_0^t \widetilde B[s, X^0_s, \mu_s^{2}]\,ds, \quad 0\le t\le T.
$$
This is a new Wiener process on $[0,T]$ under the probability measure $P^{\gamma^2}$ defined by its density as $(dP^{\gamma^2}/dP)(\omega) = \gamma^2$.
Then, on the interval $[0,T]$, on the probability space with a Wiener process $(\Omega, {\cal F}, (\widetilde W^2_t, F^{}_t), \mathbb P^{\gamma^2})$, the process $(X^0_t, \, 0\le t\le T)$ satisfies the equation,
\begin{eqnarray*}
dX^0_t
= \sigma(t, X^0_t)d \widetilde W^2_t
+B[t,X^0_t, \mu^2_t]\,dt,
\end{eqnarray*}
with the initial condition $X^0_0 =x_0$. In other words, the process $X^0$ on $[0,T]$ satisfies the equation (\ref{new2}), just with another Wiener process and under another measure. However, given $\mu^2_t, \, 0\le t\le T$, this solution considered as a solution of It\^o's equation is weakly unique. Therefore, the couple $(X^0_t, \mu^2_t)$ under the probability measure $\mathbb P^{\gamma^2}$ solves the McKean-Vlasov equation (\ref{e1}), that is, it is equivalent to the pair $(X^2_t, \mu^2_t, \, 0\le t\le T)$ under the measure $\mathbb P$.

\medskip

\noindent
{\bf 2.}
This provides us a
way to write down the density of the distribution of $X^1$ on $(\Omega, {\cal F}, \mathbb P)$ with respect to the distribution of $X^2$ on  $(\Omega, {\cal F}, \mathbb P)$ on the interval of time $[0,T]$. We have, for any measurable $A \subset C[0,T; \mathbb R^d]$,
\begin{equation}\label{mu1}
\mu_{0,T}^\mathbf{1}(A) := \mathbb P(X^1 \in A) = \mathbb P^{\gamma^1} (X^0 \in A)
= \mathbb E^{\gamma^1} \mathbf{1}(X^0 \in A)
=\mathbb E \gamma^1_T \mathbf{1}(X^0 \in A),
\end{equation}
and
\begin{equation}\label{mu2}
\mu_{0,T}^2(A) := \mathbb P(X^2 \in A) = \mathbb P^{\gamma^2} (X^0 \in A)
= \mathbb E^{\gamma^2} \mathbf{1}(X^0 \in A)
=\mathbb E \gamma^2_T \mathbf{1}(X^0 \in A).
\end{equation}
So, on the sigma-algebra ${\cal F}^W_T$ we obtain that
\begin{eqnarray*}
\displaystyle \frac{\mu^2_{[0,T]}(dX)}{\mu^1_{[0,T]}(dX)}(X^0) = \frac{\gamma^2_T}{\gamma^1_T}(X^0)
= \exp\left(\int_0^T
\widetilde B[s,X^0_s,\mu^2_s]  dW_s
- \frac12
\int_0^T |\widetilde B[s,X^0_s,\mu^2_s]|^2 ds\right)
 \\
\times
\exp\left(-\int_0^T
\widetilde B[s,X^0_s,\mu^1_s]  dW_s
+ \frac12
\int_0^T |\widetilde B[s,X^0_s,\mu^1_s]|^2 ds\right)
 \\
= \exp\left(\int_0^T
(\widetilde B[s,X^0_s,\mu^2_s] - \widetilde B[s,X^0_s,\mu^1_s]) dW_s
- \frac12
\int_0^T (|\widetilde B[s,X^0_s,\mu^2_s]|^2 - |\widetilde
B[s,X^0_s,\mu^1_s]|^2) ds\right)
 \\
=   \exp\left(\int_0^T
(\widetilde B[s,X^0_s,\mu^2_s] - \widetilde B[s,X^0_s,\mu^1_s])
\sigma^{-1}(s,X_s^0)dX^0_s\right)
\\
\times  \exp\left(- \frac12
\int_0^T (|\widetilde B[s,X^0_s,\mu^2_s]|^2 - |\widetilde
B[s,X^0_s,\mu^1_s]|^2) ds\right)
\\
= \exp\left(\int_0^T
(\widetilde B[s,X^0_s,\mu^2_s] - \widetilde B[s,X^0_s,\mu^1_s])(d\widetilde W^1_s
 + \widetilde B[s,X^0_s, \mu^1_s]\,ds) \right)
 \\
\times  \exp\left(- \frac12
\int_0^T (|\widetilde B[s,X^0_s,\mu^2_s]|^2 - |\widetilde
B[s,X^0_s,\mu^1_s]|^2) ds)\right)
\\
\displaystyle = \exp\left(\int_0^T
(\widetilde B[s,X^0_s,\mu^2_s] - \widetilde B[s,X^0_s,\mu^1_s])\,d\widetilde W^1_s
- \frac12
\int_0^T |\widetilde B[s,X^0_s,\mu^2_s] - \widetilde
B[s,X^0_s,\mu^1_s]|^2 ds\right).
\end{eqnarray*}

\medskip

\noindent
Further, due to (\ref{mu1}) and (\ref{mu2}) the measure \(\mu^i\) is an image of \(\mathbb P^{\gamma^i}\) under the mapping \(X^0\) for \(i=1,2\). So,
\begin{equation}\label{vtf}
v(t):= \|\mu^1_{[0,t]} - \mu^2_{[0,t]}\|_{TV} \le \|P^{\gamma^1}|_{{\cal F}^W_{t}} - P^{\gamma^2}|_{{\cal F}^W_{t}}\|_{TV}.
\end{equation}
Since the two measures $P^{\gamma^1}$ and $P^{\gamma^1}$ on ${\cal F}^W_{t}$ are equivalent with the density
$$
\frac{dP^{\gamma^2}}{dP^{\gamma^1}}\Big|_{{\cal F}^W_{t}} (\omega) =  \frac{\gamma^2_t}{\gamma^1_t}(\omega),
$$
the total variation distance between them equals (denoting $\rho_t = \gamma^2_t/\gamma^1_t$),
$$
\frac12\, \|P^{\gamma^2}|_{{\cal F}^W_{t}} - P^{\gamma^1}|_{{\cal F}^W_{t}}\|_{TV}
= \int_{\Omega}\left(1- \frac{\gamma^2_t}{\gamma^1_t}(\omega) \wedge 1\right) \, \mathbb P^{\gamma^1}(d\omega)
= 1 - \mathbb{E}^{\gamma^1}
\rho_t \wedge 1
 \\\nonumber\\ \nonumber
\le \sqrt{E^{\gamma^1}\rho_t^2 - 1}. \hspace{2.5cm}
$$
Let us justify the last inequality
for completeness, dropping the sub-index \(t\):
\begin{eqnarray*}
 1 - \mathbb{E}^{\gamma^1} (\rho \wedge 1) = \mathbb{E}^{\gamma^1} (1-\rho
\wedge 1) \hspace{2cm}
\\
\le  \sqrt{\mathbb{E}^{\gamma^1} (1-\rho
\wedge 1)^2} =  \sqrt{E^{\gamma^1} (1-\rho \mathbf{1}(\rho\le 1) -
\mathbf{1}(\rho>1))^2}
\\
= \sqrt{\mathbb{E}^{\gamma^1} (\mathbf{1}(\rho\le 1)-\rho \mathbf{1}(\rho\le 1))^2} =
\sqrt{\mathbb{E}^{\gamma^1} \mathbf{1}(\rho\le 1)(\rho - 1)^2}
\\
\le \sqrt{\mathbb{E}^{\gamma^1} (\rho - 1)^2} =
\sqrt{\mathbb{E}^{\gamma^1}\rho^2 - 1},  \hspace{2cm}
\end{eqnarray*}
as required. We used the CBS  inequality.
So, due to (\ref{vtf}),
\begin{eqnarray}\label{cts}
v(t)
\le 2\sqrt{\mathbb E^{\gamma^1}\rho_t^2 - 1}.
\end{eqnarray}
Now, again by virtue of the CBS  inequality,
\begin{eqnarray}\label{rolin}
\mathbb{E}^{\gamma^1} \rho_T^2 =
\mathbb{E}^{\gamma^1}
\exp\Big(- 2\int_0^T
(\widetilde B[s,X^0_s,\mu^2_s] - \widetilde B[s,X^0_s,\mu^1_s]) d\widetilde W^1_s
 \nonumber \\\nonumber \\
- \int_0^T |\widetilde B[s,X^0_s,\mu^2_s] - \widetilde
B[s,X^0_s,\mu^1_s]|^2 ds\Big)
 \nonumber \\\nonumber \\
=
\mathbb{E}^{\gamma^1}
\exp(- 2\int_0^T
(\widetilde B[s,X^0_s,\mu^2_s] - \widetilde B[s,X^0_s,\mu^1_s]) d\widetilde W^1_s
 \nonumber \\\nonumber \\
- 4 \int_0^T |\widetilde B[s,X^0_s,\mu^2_s] - \widetilde
B[s,X^0_s,\mu^1_s]|^2 ds)
 \nonumber \\\nonumber \\
\times \exp\Big(  3 \int_0^T |\widetilde B[s,X^0_s,\mu^2_s] - \widetilde
B[s,X^0_s,\mu^1_s]|^2 ds\Big)
 \nonumber \\\nonumber \\
\le \left(\mathbb{E}^{\gamma^1}
\exp\Big(- 4\int_0^T
(\widetilde B[s,X^0_s,\mu^2_s] - \widetilde B[s,X^0_s,\mu^1_s]) d\widetilde W^1_s \right.
 \nonumber \\\nonumber \\
\left.- 8 \int_0^T |\widetilde B[s,X^0_s,\mu^2_s] - \widetilde
B[s,X^0_s,\mu^1_s]|^2 ds\Big)\right)^{1/2}
 \nonumber \\\nonumber \\
\times \left(\mathbb{E}^{\gamma^1}\exp\left(
6\int_0^T |\widetilde B[s,X^0_s,\mu^2_s] -
\widetilde B[s,X^0_s,\mu^1_s]|^2 ds\right)\right)^{1/2}
 \nonumber \\\nonumber \\
\le (=) \; \sqrt{\mathbb{E}^{\gamma^1}\exp\left(
6\int_0^T |\widetilde B[s,X^0_s,\mu^2_s] -
\widetilde B[s,X^0_s,\mu^1_s]|^2 ds\right)}.
\end{eqnarray}
  The last    inequality  is always true; for a   bounded  $\widetilde B$ it is, apparently, an equality.

\medskip

\noindent
{\bf 3.}
If   $\widetilde B$ is    bounded with the norm $\|\widetilde B\|_B:=\sup_{s,\mu}|\widetilde B[s,x,\mu]|$, then
\begin{equation}\begin{gathered}\label{nest1}
 \mathbb{E}^{\gamma^1}\exp\left(
6\int_0^T |\widetilde B[s,X^0_s,\mu^2_s] - \widetilde
B[s,X^0_s,\mu^1_s]|^2 ds\right)
\\
\le \mathbb{E}^{\gamma^1}\exp\left(6 \|\widetilde B\|_B^2
\int_0^T  \|\mu^1_s -
\mu^2_s\|_{TV}^2 \,ds\right).
\end{gathered}\end{equation}
Here the value under the expectation is non-random; hence, the symbol of this expectation may be dropped. Therefore, we have with $C=6\|\widetilde B\|^2_B$,
\begin{equation}\label{neqv}
v(T) \le 2\sqrt{\exp\left(C \int_0^T v(s)^2 ds\right) -1}.
\end{equation}
Recall that $v(t)\le 2$, and the function $v$ increases in
$t$. Let us choose $\alpha_0>0$ small so that for any $0\le \alpha \le \alpha_0$,
\begin{equation}\label{al}
\exp(\alpha)-1\le 2\alpha,
\end{equation}
and take $T\le \alpha_0/(4C)$.
Then $\displaystyle C \int_0^T v(s)^2 \,ds \le CTv(T)^2 \le 4 CT \le \alpha_0$. So,
\begin{equation*}
v(T) \le 2\sqrt{\exp\left(C \int_0^T v(s)^2 \,ds\right) -1} \le 2\sqrt{2C T v(T)^2} = 2\sqrt{2CT} v(T).
\end{equation*}
If we choose $T$ so small that $2\sqrt{2CT}  < 1$, that is, $T<1/(8C)$, then it follows that $v(T) = 0$. Hence,  $v(T) = 0$ for any $T < \min(1/(8C), \alpha_0/(4C))$. Let us fix some \(T>0\) satisfying this inequality.

~

Further, we conclude by induction that
\begin{equation}\label{ind}
v(2T) = v(3T) = \ldots = 0.
\end{equation}
Indeed, assume that $v(kT)=0$ is already established for some integer $k>0$. Redefine
the stochastic exponents:
$$
\gamma^i_{kT, (k+1)T}= \exp\left( \int_{kT}^{(k+1)T} \widetilde B[s,X^0_s,\mu^i_s]\,dW_s -
\frac12
\int_{kT}^{(k+1)T} |\widetilde B[s,X^0_s,\mu^i_s]|^2 \,ds\right), \quad i=1,2,
$$
and re-denote
$$
\widetilde W^1_t := W_t  - \int_{kT \wedge t}^{t} \widetilde B[s, X^0_s, \mu_s^{1}]\,ds, \quad 0\le t\le (k+1)T.
$$
Then \(\widetilde W^1_t\) is a new Wiener process on $[kT,(k+1)T]$ starting at $W_{kT}$ under the probability measure $P^{\gamma^1}$ defined by its density as $(dP^{\gamma^1}/dP)(\omega) = \gamma^1_{kT, (k+1)T}$. Repeating the calculus leading to \eqref{rolin},  \eqref{nest1}   and  \eqref{neqv}, and having in mind the induction assumption $v(kT)=0$, we obtain  with the same constant $C$ that
\begin{equation*}
v((k+1)T) \le \sqrt{\exp\left(C \int_{kT}^{(k+1)T} v(s)^2 ds\right) -1},
\end{equation*}
which straightforwardly  implies that
\begin{equation*}
v((k+1)T) \le \sqrt{2C T v((k+1)T)^2} = \sqrt{2CT} v((k+1)T).
\end{equation*}
As earlier, the condition $T < \min(1/(2C), 1/(\alpha C))$ (see (\ref{al})) guarantees that
$$
v((k+1)T) = 0,
$$
as required. This completes the induction (\ref{ind}).

\medskip

Hence, solution is weakly unique on the whole $\mathbb R_+$. As noticed above, strong uniqueness also follows.
For bounded $\widetilde B$ the statements of Theorem  \ref{thm5} as well as of the Remark \ref{rem3} are justified.

\ifleft

\medskip

\noindent
{\bf 4.} Now let us return to the inequality (\ref{rolin}) and explain how to drop the additional assumption of boundedness of $\widetilde B$, and also how to deal with a localized version of  (\ref{liptvr}).
First of all, prior to (\ref{rolin}) we have to show that \(\gamma^{i}, \, i=1,2\), are, indeed, probability densities for which it suffices to show uniform integrability for \(T>0\) small enough: for example, it suffices to check that
\[\mathbb E (\gamma^{i}_T)^2 < \infty, \quad i=1,2.
\]
Via the estimates similar to (\ref{rolin}) by virtue of CBS inequality, this problem is reduced to the question whether or not the following expression is finite:
\begin{align}\label{ui1}
\mathbb E (\gamma^i_T)^2
\le \left(\mathbb E \exp\left(4\int_0^T \widetilde B[s,X^0_s,\mu^i_s]\,dW_s -
8 \int_0^T |\widetilde B[s,X^0_s,\mu^i_s]|^2 \,ds\right)\right)^{1/2}
 \nonumber \\
\times \left(\mathbb E \exp\left(6\int_0^T |\widetilde B[s, X^0_s, \mu^i_s]|^2\,ds\right)\right)^{1/2}
 \nonumber \\
\le  \left(\mathbb E \exp\left(6\int_0^T |\widetilde B[s, X^0_s, \mu^i_s]|^2\,ds\right)\right)^{1/2}
 \nonumber \\
\le  \left(\mathbb E \exp\left(C\int_0^T (1+|X^0_s|^2)\,ds\right)\right)^{1/2}\,.
\end{align}
In the last inequality the assumption on the linear growth of $\widetilde B$ was used.

\medskip

Suppose for instant that the finiteness of the last expectation in the last line of (\ref{ui1}) has been shown; then, by standard induction arguments with conditional expectations it follows that both \(\gamma^i_T\) are, indeed, probability densities for {\em any} \(T>0\). Hence, the calculus leading to (\ref{cts}) and (\ref{rolin}) is valid and we have that
\begin{equation}\label{nado0}
v(t) = \|\mu^1_{[0,t]} - \mu^2_{[0,t]}\|_{TV} \le 2\sqrt{\mathbb E^{\gamma^1}\rho^2 - 1},
\end{equation}
and
\begin{equation}\label{nado1}
\mathbb{E}^{\gamma^1} \rho^2
\le  \sqrt{\mathbb{E}^{\gamma^1}\exp\left(
6\int_0^T \left|\widetilde B[s,X^0_s,\mu^2_s] -
\widetilde B[s,X^0_s,\mu^1_s]\right|^2 ds\right)}.
\end{equation}
It is a general fact which does not use any boundedness of $\widetilde{B}$ in any variable but only in the last variable is,
\begin{equation}\label{blingrow}
\left|\widetilde B[s,X^0_s,\mu^2_s] -
\widetilde B[s,X^0_s,\mu^1_s]\right| \le \sup_{y}\left|\widetilde B(s,X^0_s,y)\right| \|\mu^2_s-\mu^1_s\|_{TV}.
\end{equation}
Due to the linear growth assumption (\ref{bgrow}), the inequality   (\ref{blingrow}) implies
\begin{equation}\label{blingrow2}
\left|\widetilde B[s,X^0_s,\mu^2_s] -
\widetilde B[s,X^0_s,\mu^1_s]\right| \le C(1+|X^0_s|) \|\mu^1_s-\mu^2_s\|_{TV}.
\end{equation}
Hence, by virtue of (\ref{blingrow2}) we obtain
\begin{eqnarray}\label{nado}
\mathbb{E}^{\gamma^1} \rho^2
\le  \mathbb{E}^{\gamma^1}\exp\left(
6\int_0^T [C(1+|X^0_s|) \|\mu^1_s-\mu^2_s\|_{TV}]^2 ds\right)
 \nonumber\\
\le  \mathbb{E}^{\gamma^1}\exp\left(
6C^2 v(T)^2 \int_{0}^{T}(1+|X^0_s|^2)\,ds \right).
\end{eqnarray}
Recall that the process $X^0$ satisfies the equation (\ref{x01}) on $[0,T]$ with respect to the measure $\mathbb P^{\gamma^1}$.
We want to show that given $C$, the right hand side in (\ref{ui1})
is finite for any $T$ small enough. For this end, denote  $6C^2 v(T)^2 : = r\ge 0$. Recall that in any case $v(T)\le 2$. We would like to show that for any fixed constant $K>0$ ($K = 24 C^2+1$ suffices), the value
\[
\mathbb{E}^{\gamma^1}\exp\left(
r \int_{0}^{T}(1+|X^0_s|^2)\,ds  \right)
\]
is finite for $0\le r < K$, and differentiable with respect to $r$, and that this derivative is non-negative and small uniformly in $r \in [0, K)$, if $T>0$ is  small enough.

\medskip

It suffices to show the same properties, still for small enough $T$, for the function
\begin{equation}\label{psi}
\psi(r,T) = \mathbb{E}^{}\exp\left(
r \int_{0}^{T}(1+|X^1_s|^2)\,ds \right),
\end{equation}
where $X^1$ solves the equation (\ref{x01}) on $[0,T]$ with respect to the original measure $\mathbb P^{}$, because $X^1$ solves the same equation with respect to the measure $\mathbb P$ as the process $X^0$ with respect to the measure $\mathbb P^{\gamma^1}$ on $[0,T]$.

Denote \(\bar X^{1}_t := \sup_{0\le s\le t} |X^{1}_s|\). From the equation (\ref{new1}) and due to the assumption (\ref{bgrow}) and boundedness of $\sigma$, we have that

\[
|X^{1}_t| \le |x_0| + C\left(\int_0^t (1+|X^{1}_s|)\,ds + \left|\int_0^t \sigma(s, X^1_s)\,dW^1_s\right|\right).
\]
Now, by virtue of Gronwall's inequality and since both sides in this inequality are finite, we obtain with some $C>0$,
\[
\bar X^{1}_t \le C\exp(Ct)\left(|x_0| + \sup_{0\le t'\le t} \left|\int_0^{t'} \sigma(s, X^1_s)\,dW^1_s\right|\right).
\]
Therefore, for each $r>0$
\[
\exp(r T(\bar X^{1}_T)^2)  \le
C_T\exp(2rT |x_0|^2)\exp\left(2rT \sup_{0\le t'\le T}\left|\int_0^{t'} \sigma(s, X^1_s)\,dW^1_s\right|^2\right),
\]
with some $C_T<\infty$.
Recall that the matrix-function $\sigma$ is bounded. Denote by
$
\sigma^{i}(r,x)
$
the $i$th row of the matrix $\sigma$, $1\le i\le d$.
Due to the exponential martingale inequalities following from Girsanov  theorem along with Doob  inequality for martingales (see, e.g., \cite[Sections 6.8 \& 3.4]{Kry-ln}) there exist constants $C_1, C_2>0$ such that for any $a>0$, and for any $t>0$, and for any $i$,
\[
\mathbb P(\sup_{0\le s\le t}\left|\int_0^s \sigma^i(r, X^1_r)dW^1_r\right|>a)
\le C_1\exp(-a^2/(C_2t)).
\]
This implies that (with new $C_1, C_2>0$)
\begin{align*}
\mathbb P\left(\sup_{0\le s\le T}\left|\int_0^s \sigma(r, X^1_r)dW^1_r\right|>a\right)
 \\
\le \sum_{i=1}^d \mathbb P\left(\sup_{0\le s\le T}\left|\int_0^s \sigma^i(r, X^1_r)dW^1_r\right|>\frac{a}{d}\right)
\le C_1\exp(-a^2/(C_2T)).
\end{align*}
By integration, it follows that for $2rT<1/(C_2 T)$ with $C_2$ from the last line,
\begin{align*}
\mathbb E \exp\left(2rT \sup_{0\le s\le T}\left|\int_0^s \sigma(r,X^1_r)dW^1_r\right|^2\right)<\infty,
 \end{align*}
and  hence we also have that
\begin{equation*}
\psi(r,T)\le \exp(rT)\,C_T\mathbb E \exp(2rT(1 + |x_0|^2))\,
\mathbb E \exp(2rT \sup_{0\le s\le T}|X^1_s|^2) < \infty,
\end{equation*}
as required, in particular, due to the exponential moment assumed for $x_0$.

\medskip

Thus, the function $\psi(r,T)$ (see (\ref{psi})) is finite for $r$ from some (any) finite range $0\le r<K$, if $T$ is small enough.
It follows that all expressions in (\ref{ui1}) for small enough \(T>0\) are finite. So, in particular, both \(\gamma^i_T\) are, indeed, probability densities  for small \(T>0\) under the linear growth condition (\ref{bgrow}), too. Hence, we can return to the inequalities (\ref{cts}) earlier established for bounded \(\widetilde B\), and by virtue of (\ref{nado0}), (\ref{nado1}), and (\ref{nado}) we get,
\begin{eqnarray*}
v(T) \le 2\sqrt{E^{\gamma^1}\rho^2_T - 1} \le
\sqrt{C T v(T)^2},
\end{eqnarray*}
with some constant $C$ which constant  may  depend on the initial  distribution (or value). Therefore, \(v(T)=0\) for \(T>0\) small enough.

\medskip

\noindent
{\bf 7.}
Denote
\[
{\cal N}:= \{t\ge 0: \, v(t)=0 \}.
\]
The previous steps show that \(\sup({\cal N}) > 0\) and that \(0 \in {\cal N}\). Note that \(t\in {\cal N} \Longrightarrow s\in {\cal N}\) for any \(0\le s\le t\).
Recall that \(v(t) \le \sqrt{\mathbb E^{\gamma^1}\rho_t^2 - 1}\) (see (\ref{cts})) where the right hand side is clearly continuous in \(t\). Moreover, as it follows from (\ref{nado0}), (\ref{nado1}), and (\ref{nado}),
\[
v(t)^2 \le
4(\mathbb E^{\gamma^1}\rho_t^2 - 1)
\le  4\left(\mathbb{E}^{\gamma^1}\exp\left(
6\int_0^t [C(1+|X^0_s|) \|\mu^1_s-\mu^2_s\|_{TV}]^2 ds\right) -1\right),
\]
which implies that the set \({\cal N}\) is closed, since $\|\mu^1_s-\mu^2_s\|_{TV}=0$ for all $s<t$ if $[0,t)\subset {\cal N}$.

\medskip

On the other hand, consider any  \(N \in (0, \sup({\cal N}))\). Recall that \(\mathbb E \exp(c \,\sup_{s\le N} |X^0_s|^2) < \infty\) with some positive \(c\).
Hence, the same calculus as above shows that \(v(t)=0\) in some small right neighbourhood of \(N\). In other words, the set on the positive half-line \(\mathbb R_+\) where \(v(t)=0\) is non-empty, closed and open in \(\mathbb R_+\). Thus, it coincides with \(\mathbb R_+\) itself. In other words, for all  \(t\ge 0\),
\[
v(t)=0,
\]
which finishes the proof of Theorem  \ref{thm5a}.
\fi

\subsection{Strong uniqueness}
In this section it will be shown that in certain cases weak uniqueness
implies strong uniqueness for the equation (\ref{e1}) -- (\ref{e200}), and  both properties will be established under appropriate
conditions. This result, Theorem \ref{thm5} below, requires only a Borel  measurability
of the drift with respect to the state variable $x$, although, it
assumes that diffusion $\sigma$ does not depend on $y$
along with Lipschitz condition in $x$ and non-degeneracy.
The drift may be unbounded in the state variable $x$.

\begin{theorem}\label{thm5}
	Let
	\(\mathbb E \exp(r |x_0|^2)<\infty\) for some \(r>0\), and let the functions
	$b$ and $\sigma$ be Borel
	measurable,  and
	$$
	\sigma(s,x,y) \equiv \sigma (s,x),
	$$
	that is, \(\sigma\) does not depend on the variable \(y\);
	let $\sigma$ satisfy the non-degeneracy assumption (\ref{si1});
	let $d_1 = d$, the matrix $\sigma$ be bounded, symmetric and invertible, and let there exist $C>0$ such that the function
	$$
	\widetilde
	B[s,x,\mu] := \sigma^{-1}(s,x)\, B[s,x,\mu]
	$$
	satisfies the linear growth condition: there is $C>0$ such that for all $x\in R^d$,
	\begin{equation*}
	\sup_{s,\mu}|\widetilde B[s,x,\mu]| \le C(1+|x|).
	\end{equation*}
	Also assume that the matrix-function $\sigma(t,x)$ satisfies the following Lipschitz condition (for simplicity) which guarantees that the equation
	\begin{equation*}
	dX^0_t = \sigma(t,X^0_t)\,dW_t, \quad X^0_0 = x_0,
	\end{equation*}
	has a unique strong solution for any $x$ (see  \cite{Ver80, Ver82}): \\
	\begin{equation}\label{Lip1}
	\sup_{t\ge 0}\, \sup_{x,x': \, x'\not = x}\, \frac{\|\sigma(t,x)-\sigma(t,x')\|}{|x-x'|}<\infty.
	\end{equation}
	
	Then  solution of the equation (\ref{e1})--(\ref{e200}) is weakly and
	strongly unique; this solution is strong.
\end{theorem}

\noindent
\begin{proof} It follows straightforwardly from the Theorem \ref{thm5a} and from the fact that with a given $\mu_t$ any solution is strong \cite{Ver80, Ver82} (note that linear growth of the drift is allowed in both \cite{Ver80, Ver82}).
\end{proof}

\begin{remark}\label{rem2b}
	Note that under the condition (\ref{Lip1}), not only the equation (\ref{itosde}) but any equation with the same diffusion and a Borel measurable drift with a linear growth assumption in $x$ will have a strong solution. It concerns both solutions of the equation (\ref{e1}) and its ``linearized'' version (\ref{e22}).

	Emphasize that no regularity on the function $b$ is needed in either variable. Also, a linear growth condition on the drift in $x$ is equivalent  to the condition (\ref{bgrow}); the latter was assumed in order to make the presentation more explicit. The price for the no regularity and linear growth is a special form of $\sigma$ which may not depend on
	the ``measure variable'' $y$; in particular, in such a case $\Sigma(t,x) =\sigma(t,x)$, and we will use the lower case to denote the diffusion coefficient in the remaining sections.
\end{remark}

\begin{remark}\label{rem7}
	Instead of Lipschitz condition (\ref{Lip1}), it suffices if diffusion coefficient \(\sigma\)  belongs to the Sobolev class \(\sigma(t,x) \in W^{0,1}_{2d+2, loc}\). More general conditions on Sobolev derivatives for \(\sigma\) can be found in \cite[Theorem 1]{Ver80} and \cite{Ver82}, and any of them can be used in our Theorem \ref{thm5} above. Note that in the latter reference \(\sigma\) is assumed Lipschitz but it is shown that continuity is necessary only with respect to the state variable \(x\), which is also applied to the conditions from \cite{Ver80}. As usual, a more relaxed conditions on $\sigma$ can be stated in the case of dimension one as in \cite{YamadaWatanabe}, or in \cite[Theorem~2]{Ver80}, the simplest version of both being just H\"older 1/2.
\end{remark}

\section{Appendix}\label{Appe1}
To prove Theorem \ref{thm1}, we need  two auxiliary lemmas. Lemma \ref{lebds}, is about  very standard  a priori moment bounds for solutions; another one,  Lemma \ref{lekrybd}, is a localized version of Krylov's bounds where coefficients are locally bounded and the obtained upper bound  also relates to a bounded domain.

\begin{lemma}\label{lebds}
Under   assumption (\ref{linear}) and the standard measurability,   a priori estimates
\begin{align}
 \label{kol11}
\mathbb{E} \sup_{0\le t\le T} |X_t^{}|^4 \le C_T (1+\mathbb E |x_0|^4),
\end{align}
and
\begin{eqnarray}\label{kol20}
\sup_{0\le s\le t\le T; \, t-s\le h}\mathbb{E} |X_t^{} - X^{}_s|^4 \le
C_{T}  h^2,
\end{eqnarray}
 hold  true with some constants $C_{T}$ (generally, different) that do not depend on $n$. In   (\ref{kol11}) the constant $C_T$, generally speaking, also depends on  the value of the moment $\mathbb E |x_0|^4$.
\end{lemma}

\begin{remark} In the proof of the Proposition  \ref{pro22}   similar a priori bounds  are stated for the successive approximations of solutions.
Recall that $x_0\in \mathbb R^d$ is the initial value of the process $X$ and that it may be random with a certain finite moment.
In fact, similar a priori bounds hold true for any power function assuming the appropriate initial moment, although, this will not be used in this paper. The proof of   (\ref{kol11})
 is very standard and can be done following the lines in  \cite[Theorem 1.6.4]{GS68}, or \cite[Corollary 2.5.6]{Kry}, or \cite{Sko} (and many other places) combined with  Doob's inequalities. The bounds   (\ref{kol20}) follow from a similar calculus starting from $s$ instead of $0$.

Why do we need the 4th moment,  is clarified in the proofs of the Proposition  \ref{pro22}  and Theorem \ref{thm1}: it is useful for verifying continuity for the processes with equivalent finite-dimensional distributions; for this purpose the 2nd  moment is not sufficient, although $2+\epsilon$ should probably work.
\end{remark}
The next lemma is a version of Krylov's bounds.
Let $Z_t$ be a non-explosive strong Markov process in $\mathbb R^d$
 satisfying an SDE
$$
dZ_t = b_t(Z_t)dt + \sigma_t(Z_t) dW_t, \quad Z_0 = z_0,
$$
where (non-random) functions $b_t(z)$ and $\sigma_t(z)$ are $d$-vector and matrix $d\times d$ respectively, $\sigma_t$ is uniformly non-degenerate, and locally in $z$ bounded uniformly in $t$, that is,
\begin{equation}\label{locbd}
d(R):=\sup_{|z|\le R} \sup_{t}(|b_t(z)| + \|\sigma_t(z)\|)<\infty, \quad \forall \;\; R>0,
\end{equation}
and the random variable $z_0$ is independent of the filtration of the Wiener process $W=(W_t,{\cal F}_t)$.
Let $D$ be a bounded domain in  $\mathbb R^d$,
 $g: \mathbb R^d \to \mathbb R^1$, and $f: \mathbb R^+\times \mathbb R^d \to \mathbb R^1$ be Borel measurable functions.
Then we state the   versions of Krylov's bounds in which  only local boundedness of the coefficients is assumed, and the statements are also local.

\begin{lemma}\label{lekrybd} Let $D \subset B_{R}=\{z\in R^d: |z|\le R\}$.
Under assumption (\ref{locbd}), for any $p\ge d$ there exists a constant $N=N_R$ which also depends on $d$, on the constants of ellipticity of $\sigma\sigma^*$ and on the upper bounds for the norms of $b$ and $\sigma$ for $|x|, |y| \le R$, such that for any $g$ satisfying $g(z)\equiv 0, \, z \not \in D$,
\begin{equation}\label{krhomo}
\mathbb E \int_0^{T} g(Z_t)dt \le N_R \|g\|_{L_p(D)},
\end{equation}
 and for any $f$ satisfying $f(t,z)\equiv 0, \,  z \not \in D$,
\begin{equation}\label{krinhomo}
\mathbb E \int_0^{T} f(t,Z_t)dt \le N_R \|f\|_{L_{p+1}([0,T]\times D)}.
\end{equation}
 \end{lemma}
Notice, that in  the proof of  Theorem \ref{thm1} the role of $Z$ was played by the pairs $(\widetilde X^n, \widetilde \xi^n)$. These are yet not full Krylov's estimates for general It\^o processes; however, they suffice for our goals in the present paper.

\medskip

\begin{proof} The proof is based on Krylov's bounds for It\^o processes with bounded coefficients and from the hint similar to the one in the proof of the Lemma 4.3.1 \cite{Khasminsky}.
We only establish  bound  \eqref{krinhomo} and  the bound (\ref{krhomo}) follows similarly.
Let $D'$ be another bounded domain containing the closure of $D$: $  \overline{D} \subset  D'$. Without loss of generality, we may assume that $\inf_{y\in\overline{ D}, z\in \overline{\partial D'}}|y-z|>0 $  and  $D' \subset B_{R+1}=\{z\in R^d: |z|\le R+1\}$.
Denote
$$
\tau^0=0, \quad T^1:= \inf(t\ge \tau_0: Z_t \not \in \bar D'), \quad
$$
$$
\tau^{k}:= \inf(t\ge T^k: Z_t \in D), \;
T^{k+1}:= \inf(t\ge \tau^{k}: Z_t \not \in \bar D'), \; k\ge 1,
$$
and let
\begin{align*}
\widehat Z^0_t:= Z_0 + \int_0^t \mathbf {1}(s<T^1) b_s(Z_s) ds + \int_0^t (\mathbf {1}(s<T^1)\sigma_s(Z_s)
 + \mathbf {1}(s\ge T^1))dW_s,
\end{align*}
for $t\ge 0$, and by induction,
\begin{align*}
\widehat Z_t^k:= Z_{\tau^k} + \int_{\tau^k}^t \mathbf {1}(s<T^{k+1}) b_s(Z_s) ds + \int_{\tau^k}^t (\mathbf {1}(s<T^{k+1})\sigma_t(Z_s) + \mathbf {1}(s\ge T^{k+1}))dW_s,
\end{align*}
for $t\ge \tau^k$, $k\ge 0$.
(That is, each process $\widehat Z^k$ starts in $\tau^k$ at state $Z_{\tau^k}$ and follows the trajectory of $Z$ until $T^{k+1}$, after which  continues as a Wiener process.)
Rigorously, recall that $Z_t = \widehat Z^0_t$ on the set $(t<T^1)$, and $\widehat Z_t^k = Z_t$ on $(\tau^k \le t<T^{k+1})$
(see  \cite[Theorem 6.3.7(iii)]{Kry-ln}).
Since $f(s,Z_s) = 0$ on any interval  $T^k\le s \le \tau^{k}$, we have, according to Krylov's bound \cite[Theorems 2.2.2, 2.2.4]{Kry}, that for  $p\ge d$
\begin{align*}
\mathbb E
\int_0^T |f(s,Z_s)|\,ds
= \sum_{k=0}^{\infty} \mathbb E \int_{\tau^k \wedge T}^{T^{k+1} \wedge T}\!\!\!\! |f(s,Z_s)|\,ds
= \sum_{k=0}^{\infty} \mathbb E \mathbf {1}(\tau^k \le T) \int_{\tau^k \wedge T}^{T^{k+1} \wedge T} \!\!\!|f(s,Z_s)|\,ds
 \\
= \sum_{k=0}^{\infty} \mathbb E \mathbf {1}(\tau^k \le T) \int_{\tau^k}^{T^{k+1} \wedge T} |f(s,Z_s)|\,ds
\le \sum_{k=0}^{\infty} \mathbb E \mathbf {1}(\tau^k \le T) \int_{\tau^k }^{T^{k+1}} |f(s,Z_s)|\,ds
 \\
= \sum_{k=0}^{\infty} \mathbb E \mathbf {1}(\tau^k \le T) \mathbb E \left(\int_{\tau^k }^{T^{k+1}} \!\!\!\!\!|f(s,\widehat Z^k_s)|\,ds |{\cal F}_{\tau^k }\right)
\le  N_R \|f\|_{L_{p+1}([0,T]\times R^d)} \sum_{k=0}^{\infty} \mathbb{P}(\tau^k \le T).
\end{align*}
 Recall that $\mathbb{P}(\tau^k \le T) \le \mathbb{P}(T^k \le T)$. For a strong Markov process $Z$ with a positive probability to exit from $\bar D'$ on any finite interval of time due to the non-degeneracy of its diffusion coefficient and boundedness of both coefficients in $\bar D'$, the probabilities $\mathbb{P}(T^k \le T)$ admit exponential bounds
 $$
\mathbb{P}(T^k \le T) \le C q^{k-1}, \; k\ge 1,
$$
with some $C<\infty$ and  $q<1$.
So, with a new constant $C$ and since $( {1}-\mathbf {1}(D)) f \equiv 0$,
 \begin{align*}
\mathbb E
\int_0^T |f(s, Z_s)|\,ds
\le N_R \|f\|_{L_{p+1}([0,T]\times  \mathbb R^d)} \sum_{k=0}^{\infty}   \mathbb{P}(T^k \le T)
  \le C N_R\|f\|_{L_{p+1}([0,T]\times D)},
\end{align*}
as required. \hfill \end{proof}

\begin{lemma}[Skorokhod \mbox{(on unique probability space and convergence)}] \label{app2}

Let $\{\xi^n_t,\;t\ge 0,\;n\ge 0\}$ be some $d$-dimensional  stochastic processes defined on some probability space (the spaces can be different for different $n$) and let for any $T>0$, $\varepsilon > 0$ the following hold true:
\[
\lim_{c\to\infty} \sup_n \sup_{t\le T}\mathbb P(|\xi^n_t|>c) = 0,
\]
and
\[
\lim_{h\downarrow \infty} \sup_n \sup_{t,s\le T; \, |t-s|\le h}\mathbb P(|\xi^n_t - \xi^n_s|>\varepsilon) = 0,
\]
Then for any sequence $n'\rightarrow \infty$   a new probability space $(\Omega',\mathcal{F}', \mathbb{P}')$can be constructed that supports  the  processes $\widetilde \xi^{n'}_t, \, t\ge 0$ and $\widetilde \xi^{}_t, \, t\ge 0$, such that all finite-dimensional distributions of $\widetilde \xi^{n'}_{\cdot}$ coincide with those of $\xi^{n'}_{\cdot}$ and there exists a subsequence $n\to\infty$ such that  for any $t\ge 0$
\[
 \widetilde \xi^{n}_t \rightarrow \widetilde \xi^{}_t,  \quad n\to\infty,
\]
in probability ${\mathbb{P}}'$.
\end{lemma}
See
\cite[Ch.1, \S 6]{Sko}.

\begin{lemma}[Skorokhod] \label{app1}
Let  $f^n: \mathbb R\times\Omega\rightarrow \mathbb R, \,  n\geq 0 $ be uniformly bounded random processes on some probability space; let $(W^n, \,  n\ge 0 )$ be a sequence of  one-dimensional  Wiener processes on the same probability space, and let all It\^o's stochastic integrals $\displaystyle \int_0^T f^n_s\, dW^n_s, n\geq 0$   be well-defined. Assume that for any $\varepsilon>0$,
\begin{equation*}
\lim_{h\rightarrow 0}\sup_{n}\sup_{|s-t|\leq h}\mathbb{P}\{|f^n_s-f^n_t|>\varepsilon\}=0,
\end{equation*}
and let for each $s\in [0,T]$
$$
(f^n_s, W^n_s)\stackrel{\mathbb{P}}{\to} (f^0_s, W^0_s).
$$
 Then
$$
\int_0^T f^n_s\, dW^n_s  \stackrel{\mathbb{P}}{\to} \int_0^T f^0_s\, dW^0_s.
$$
\end{lemma}
See \cite[Ch.2, \S 3, Theorem]{Sko}, where \(W^n\) are allowed to be more general martingales with brackets converging to that of a Wiener process.

\medskip

\ifleft

\fi

\section*{Acknowledgements}
  For the second author this research  has been funded
by the Russian Academic Excellence Project '5-100' (Proposition \ref{pro22}, Lemma \ref{lekrybd}) and by the Russian Science Foundation project 17-11-01098 (Theorem \ref{thm5a}).
Certain stages of this work have been fulfilled while the second author was visiting Bielefeld university to which programme SFB1283 this author is grateful.

Both authors are grateful to Denis Talay, Mireille Bossy and Sima Mehri who posed valuable questions and thus drew attention of the authors to various gaps in the earlier proofs mainly in the Theorem \ref{thm1} and in this way stimulated us to improve our presentation. Their help was indispensable. The authors are grateful to the referee for many useful comments and remarks which helped essentially with final corrections.

\end{document}